\documentclass[preprint,3p,12pt]{elsarticle}
\usepackage{mathrsfs}
\usepackage{amsmath}
\usepackage{stmaryrd}
\usepackage{bbding}
\usepackage{dcolumn}
\usepackage{graphicx}
\usepackage{amsfonts}
\usepackage{amssymb}
\usepackage{psfrag}
\usepackage{wrapfig}
\usepackage{subfigure}
\usepackage{makeidx}
\usepackage{bm}
\usepackage{epsf}
\usepackage{epsfig}
\usepackage{setspace}
\usepackage{graphicx}
\usepackage{epstopdf}
\usepackage{psfrag}
\usepackage{subfigure}
\usepackage{color}

\newtheorem{rmk}{Remark}

\begin{document}
\title{Weighted essentially non-oscillatory scheme on unstructured quadrilateral and triangular meshes for hyperbolic conservation laws}

\author[caep,iapcm]{Fengxiang Zhao}
\ead{kobezhao@126.com}

\author[iapcm]{Liang Pan\corref{cor}}
\ead{panliangjlu@sina.com}

\author[iapcm]{Shuanghu Wang}
\ead{wang\_shuanghu@iapcm.ac.cn}

\address[caep]{The Graduate School of China Academy of Engineering Physics, Beijing, China}
\address[iapcm]{Institute of Applied Physics and Computational Mathematics, Beijing, China}
\cortext[cor]{Corresponding author}

\begin{abstract}
In this paper, a third-order weighted essentially non-oscillatory
(WENO) scheme is developed for hyperbolic conservation laws on
unstructured quadrilateral and triangular meshes. As a starting
point, a general stencil is selected for the cell with any local
topology, and a unified linear scheme can be constructed. However,
in the traditional WENO scheme on unstructured meshes, the very
large and negative weights may appear for the mesh with lower
quality, and the very large weights make the WENO scheme unstable
even for the smooth tests. In the current scheme, an optimization
approach is given to deal with the very large linear weights, and
the splitting technique is considered to deal with the negative
weights obtained by the optimization approach. The non-linear weight
with a new smooth indicator is proposed as well, in which the local
mesh quality and discontinuities of solutions are taken into account
simultaneously. Numerical tests are presented to validate the
current scheme. The expected convergence rate of accuracy is
obtained, and the absolute value of error is not affected by mesh
quality. The numerical tests with strong discontinuities validate
the robustness of current WENO scheme.
\end{abstract}


\maketitle
\section{Introduction}
In recent decades, there have been continuous interests and efforts
on the development of high-order schemes. For engineering
applications, the construction of high-order numerical schemes on
unstructured meshes becomes extremely demanding due to the complex
computational domains. There are a gigantic amount of publications
about the introduction and survey of high-order schemes, including
discontinuous Galerkin (DG), spectral volume (SV), spectral
difference (SD), correction procedure using reconstruction (CPR),
essential non-oscillatory (ENO), weighted essential
non-oscillatory (WENO), Hermite WENO (HWENO) and $k$-exact method, etc.

The DG scheme was first proposed
\cite{DG1} to solve the neutron transport equation. The major
development of DG method was carried out for the compressible Euler
equations \cite{DG2,DG3,DG5}. High-order accuracy is achieved by
means of high-order polynomial approximation within each element
rather than by means of wide stencils. Because only neighboring
elements interaction is included, it becomes efficient in the
application of complex geometry. Recently, a correction procedure
via reconstruction framework (CPR) was developed
\cite{CPR1,CPR2,CPR3,CPR4}. The CPR formulation is based on a nodal
differential form, with an element-wise continuous polynomial
solution space. By choosing certain correction functions, the CPR
framework can unify several well known methods, such as the DG, SV
\cite{SV}, and SD \cite{SD} methods, and lead to simplified versions
of these methods, at least for linear equations. The ENO scheme
was proposed in \cite{Harten,ENO-Shu} and
successfully applied to solve the hyperbolic conservation laws and
other convection dominated problems. The ``smoothest" stencil is
selected among several candidates in reconstruction to achieve high
order accuracy and keep essentially non-oscillatory near
discontinuities. Based on unstructured triangular meshes, the ENO
scheme was developed as well
\cite{Harten-unstructured,Abgrall,Sonar}. Following the ENO scheme,
the WENO scheme \cite{WENO-Liu,WENO-JS,WENO-M,WENO-Z,Zhao} was developed. With the nonlinear convex combination of candidate
polynomials, WENO scheme achieves higher order of accuracy and keeps
essentially non-oscillatory property. WENO scheme improves upon ENO
scheme in robustness, smoothness of fluxes, steady-state
convergence, provable convergence properties, and more efficiency.
The classical development of WENO scheme includes WENO-JS scheme
\cite{WENO-JS}, WENO-M scheme \cite{WENO-M}, WENO-Z scheme
\cite{WENO-Z}, et al. Based on the Hermite interpolation, 
HWENO schemes were developed \cite{HWENO1,
HWENO2,HWENO3} for the hyperbolic conservation laws. 
On the unstructured meshes, the WENO scheme
with the same order of accuracy as the fundamental ENO scheme was
constructed by Friedrichs \cite{Friedrich}. Compared with the
fundamental ENO scheme, the WENO scheme is very smooth and stable in
smooth regions, while the optimal order of accuracy is not obtained
\cite{WENO-JS}. Hu and Shu \cite{Hu-Shu} presented the third-order
and fourth-order WENO schemes on unstructured triangular meshes. The
high-order of accuracy is obtained by the combination of lower order
polynomials, which is similar with the one-dimensional WENO scheme.
However, its successful application was limited by the appearance of
the negative linear weights, and the problem appears commonly on the
unstructured meshes. For a mesh that is close to the regular meshes,
such WENO scheme works well by a regrouping approach to avoid
negative weights. However, for a mesh with lower quality, the large
linear weights appears and WENO schemes become unstable even for the
smooth flows. A splitting technique dealing with the negative linear
weights was proposed in \cite{splitting-weights}. This technique
makes the WENO schemes \cite{Hu-Shu} generally work well on better
quality meshes. However, the linear weights can be very large with
the lower quality mesh, and the WENO schemes become unstable even if
the splitting technique is adopted. In order to avoid the negative
linear weights and very large linear weights, a class of WENO
schemes were proposed in \cite{Dumbser}. Due to the very large
stencils, the boundary treatment becomes complex and consuming
computer memory, and the computational efficiency is reduced because
of the inefficient combination in reconstruction. Based on the two
classes of WENO schemes in \cite{Hu-Shu,Dumbser}, a hybrid scheme is
presented in \cite{hybrid-WENO} to deal with the problem with very
large linear weights, while the scheme contains artificial parameter
and the reconstruction is more complex. Another type of higher-order
finite volume scheme is based on the so-called $k$-exact reconstruction
\cite{k-exact1,k-exact2,k-exact3}. To guarantee the non-singularity of
the reconstruction, the $k$-exact reconstruction relies on the
least-square approach.

In this paper, a third-order WENO scheme is developed on
unstructured quadrilateral and triangular meshes for the hyperbolic
conservation laws. As a starting point of WENO reconstruction, a
general stencil is selected for a cell with any local topology. With
the selected stencil, a unified linear scheme was constructed. With
the linear candidate polynomials from the candidate sub-stencils,
the linear wights can be obtained. However, the appearance of very
large linear wights on unstructured meshes causes the WENO schemes to
be unstable \cite{Hu-Shu,hybrid-WENO}. For both quadrilateral and
triangular meshes, an optimization approach is introduced to deal
with the very large linear weights. The splitting technique is
considered to deal with the negative weights obtained by the
optimization approach. The non-linear weight with a new smooth
indicator is defined for the solutions with discontinuities. With
the optimization approach for very large weights and the splitting
technique for negative weights, the scheme becomes more robust
compared with the previous WENO scheme \cite{Hu-Shu} with the low
quality meshes.  A variety of numerical tests from the accuracy test
to the solutions with strong discontinuities are presented to
validate the accuracy and robustness of the current scheme.

This paper is organized as follows. In Section 2, the finite volume
scheme on unstructured meshes are introduced. The third-order WENO
schemes on unstructured quadrilateral and triangular meshes are
presented in section 3. Section 4 includes numerical tests to
validate the current algorithm. The section 5 is the conclusion.

\section{Finite volume method}
In this paper, the two-dimensional Euler equations are considered
\begin{equation}\label{hyperbolic}
\frac{\partial W}{\partial t}+ \frac{\partial F(W)}{\partial x}+
\frac{\partial G(W)}{\partial y}=0,
\end{equation}
where $W=(\rho, \rho U, \rho V, \rho E)$, $F(W)=(\rho U, \rho U^2+p,
\rho UV, (\rho E+p)U)$, $G(W)=(\rho V, \rho UV, \rho V^2+p, (\rho
E+p)V)$, $\rho E=\displaystyle\frac{1}{2}\rho
(U^2+V^2)+\frac{p}{\gamma-1}$ and $\gamma$ is the specific heat
ratio. For a polygon cell $\Omega_i$, the boundary can be expressed as
\begin{equation*}
\displaystyle \partial \Omega_i=\bigcup_{p=1}^n\Gamma_{ip}.
\end{equation*}
where $n$ is the number of cell interface for cell $\Omega_i$.
Integrating Eq.\eqref{hyperbolic} over the cell $\Omega_i$, the
semi-discretized form of finite volume scheme can be written as
\begin{equation}\label{semidiscrete}
\frac{\text{d}W_{i}}{\text{d}t}=\mathcal{L}(W^{n}_i)=-\frac{1}{\left| \Omega_i \right|} \sum_{p=1}^n\int_{\Gamma_{ip}}
\textbf{F}(W)\cdot\textbf{n}_p \text{d}s,
\end{equation}
where $W_{i}$ is the cell averaged value over cell $\Omega_i$, $\left|
\Omega_i \right|$ is the area of $\Omega_i$, $\textbf{F}=(F,G)^T$
and $\textbf{n}_p$ is the outer normal direction of $\Gamma_{ip}$.

In this paper, a third-order spatial reconstruction will be
introduced, and the line integral over $\Gamma_{ip}$ is discretized
according to Gaussian quadrature as follows
\begin{equation}\label{quadrature}
\int_{\Gamma_{ip}}\textbf{F}(W)\cdot\textbf{n}_p
\text{d}s=\left|l_p\right| \sum_{k=1}^{2} \omega_k
\textbf{F}(\textbf{x}_{p,k},t)\cdot \textbf{n}_p,
\end{equation}
where $\left|l_p\right|$ is the length of the cell interface
$\Gamma_{ip}$, $\displaystyle\omega_1=\omega_2=1/2$ are the Gaussian
quadrature weights, and the Gaussian quadrature points
$\textbf{x}_{pk}, k=1,2$ for $\Gamma_{ip}$ are defined as
\begin{equation*}
\textbf{x}_{p1}=\frac{3+\sqrt{3}}{6}\textbf{X}_{p1}+\frac{3-\sqrt{3}}{6}\textbf{X}_{p2},~
\textbf{x}_{p2}=\frac{3+\sqrt{3}}{6}\textbf{X}_{p2}+\frac{3-\sqrt{3}}{6}\textbf{X}_{p1},
\end{equation*}
where $\textbf{X}_{p1}, \textbf{X}_{p2}$ are the endpoints of $\Gamma_{ip}$.

For the two-dimensional Euler equations, the rotational invariance
property \cite{Riemann-Solver} can be expressed as
\begin{align}\label{transform}
\textbf{F}(W)\cdot\textbf{n}=
F(W)\cos\theta+G(W)\sin\theta=T^{-1}F(TW) ,
\end{align}
where $T=T(\theta)$ is the rotation matrix
\begin{equation*}
T= \left(
  \begin{array}{cccc}
    1 & 0         & 0         & 0 \\
    0 & \cos\theta & \sin\theta & 0 \\
    0 &-\sin\theta & \cos\theta & 0 \\
    0 & 0         & 0         & 1 \\
  \end{array}
\right).
\end{equation*}
The numerical fluxes in Eq.\eqref{quadrature} at the Gaussian quadrature points
in the global coordinate can be given with the fluxes in the local
coordinate according to Eq.\eqref{transform}, and the flux
$F(\textbf{x}_{p,k},t)$ in the local coordinate can be approximated by Riemann solvers
\begin{equation*}
F(\textbf{x}_{p,k},t)=\mathbb{F}(W^+(\textbf{x}_{p,k},t),W^-(\textbf{x}_{p,k},t)).
\end{equation*}
where $W^+(\textbf{x}_{p,k},t)$ and $W^-(\textbf{x}_{p,k},t)$ are
the reconstructed variables at both sides of  Gaussian quadrature
points of cell interface. In this paper, the HLLC Riemann solver
\cite{Riemann-Solver} is adopted for numerical tests.

\section{Third-order WENO scheme on unstructured meshes}
The reconstruction procedure is the main theme of this paper, and it
will be presented on both unstructured quadrilateral and triangular
meshes in this section.

\subsection{Linear reconstruction on quadrilateral meshes}
As a starting point of WENO reconstruction, a linear reconstruction
will be presented. For
a piecewise smooth function $W(x,y)$ over cell $\Omega_{i}$, a
polynomial $P^r(x,y)$ with degree $r$ can be constructed to
approximate $W(x,y)$ as follows
\begin{equation*}
P^r(x,y)=W(x,y)+O(\Delta x^{r+1},\Delta y^{r+1}).
\end{equation*}
In order to achieve the third-order accuracy and satisfy
conservative property, the following quadratic
polynomial on the cell $\Omega_{i_0}$ is constructed
\begin{equation}\label{qua-def}
P^2(x,y)=W_{i_0}+\sum_{k=1}^5a_kp^k(x,y),
\end{equation}
where $W_{i_0}$ is the cell average value of $W(x,y)$ over cell $\Omega_{i_0}$ and
$p^k(x,y), k=1,...,5$ are basis functions, which are given as
follows
\begin{align}\label{base}
\begin{cases}
\displaystyle p^1(x,y)=x-\frac{1}{\left| \Omega_{i_0} \right|}\displaystyle\iint_{\Omega_{i_0}}x dxdy, \\
\displaystyle p^2(x,y)=y-\frac{1}{\left| \Omega_{i_0} \right|}\displaystyle\iint_{\Omega_{i_0}}y dxdy, \\
\displaystyle p^3(x,y)=x^2-\frac{1}{\left| \Omega_{i_0} \right|}\displaystyle\iint_{\Omega_{i_0}}x^2 dxdy, \\
\displaystyle p^4(x,y)=y^2-\frac{1}{\left| \Omega_{i_0} \right|}\displaystyle\iint_{\Omega_{i_0}}y^2 dxdy, \\
\displaystyle p^5(x,y)=xy-\frac{1}{\left| \Omega_{i_0}
\right|}\displaystyle\iint_{\Omega_{i_0}}xy dxdy.
\end{cases}
\end{align}

\begin{figure}[!htb]
\centering
\includegraphics[width=0.45\textwidth]{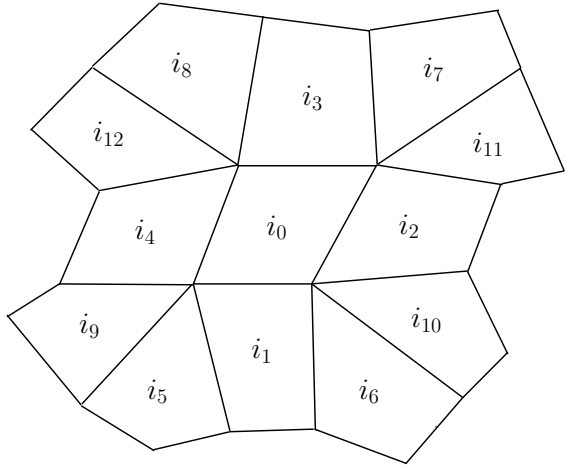}
\caption{\label{quadstencil} Stencil for third-order scheme on
quadrilateral meshes.}
\end{figure}

In order to construct the polynomial $P^2(x,y)$ uniquely, the
stencil for reconstruction is given in Figure.\ref{quadstencil}. For
arbitrary unstructured quadrilateral meshes, the cells
$\Omega_{i_j},j=1,...,4$ must exist, and the cells
$\Omega_{i_j},j=5,...,12$ may not exist since it is common that only
three or four cells share a common node. Thus, preconditioning is
needed for a general quadrilateral mesh, and the procedure of
selecting stencil for cell $\Omega_{i_0}$  is given as follows
\begin{enumerate}
\item initialize the parameter $E_j=1$, which is corresponding to the cell
$i_j, j=1,...,12$;
\item read the neighboring cell $i_1, i_2, i_3, i_4$ of cell $i_0$ from
the mesh file;
\item read the neighboring cell $i_5, i_6$ of cell $i_1$, if each of them
is not exist, record $E_j=0$ and $i_{j}$ is replaced by $i_0$,\\
read the neighboring cell $i_7, i_8$ of cell $i_3$, if each of them
is not exist, record $E_j=0$ and $i_{j}$ is replaced by $i_0$;
\item read the neighboring cell $i_{10}, i_{11}$ of cell $i_2$, if each of
them is not exist, record $E_j=0$ and $i_{j}$  is replaced by
$i_{j-4}$,\\ read the neighboring cell $i_9, i_{12}$ of cell $i_4$,
if each of them is not exist, record $E_j=0$ and $i_{j}$  is
replaced by $i_{j-4}$.
\end{enumerate}
To determine the
polynomial $P^2(x,y)$, the following constrains need to be satisfied
for the cells $\Omega_{i_j}, j=1,...,12,$
\begin{align}\label{reconstruction}
\frac{1}{\left| \Omega_{i_j}
\right|}\int_{\Omega_{i_j}}P^2(x,y)\text{d}x\text{d}y=W_{i_j},
\end{align}
where $W_{i_j}$ is the cell averaged value over $\Omega_{i_j}$.
According to the definition of $P^2(x,y)$, Eq.\eqref{reconstruction}
can be written into the following linear system
\begin{equation*}
\textbf{A}\cdot \textbf{a}=\Delta\textbf{W},
\end{equation*}
where $\textbf{a}=(a_1,...,a_5)^T$, $\Delta\textbf{W}=(W_{i_1}-W_{i_0},...,W_{i_{12}}-W_{i_0})^T$,
and $\textbf{A}=(A_{jk}),~j=1,...,12,~k=1,...,5$ is given as
\begin{equation*}
A_{jk}=E_j\cdot \frac{1}{\left| \Omega_{i_j}
\right|}\displaystyle\iint_{\Omega_{i_j}}p^k(x,y)\text{d}x\text{d}y.
\end{equation*}
Traditionally, the following least square method is used to obtain the coefficient
\begin{equation}\label{quadratic-1}
\textbf{A}^T\textbf{A}\cdot \textbf{a}=\textbf{A}^T\Delta\textbf{W}.
\end{equation}
However, the matrix $\textbf{A}^T\textbf{A}$ may not be invertible.
In order to obtain the coefficients $a_k$ in the quadratic
polynomial uniquely, a more robust least square method is adopted
\begin{equation} \label{quadratic-2}
(\zeta \textbf{I}+\textbf{A}^T\textbf{A})\cdot\textbf{a}=\textbf{A}^T\Delta\textbf{W},
\end{equation}
where $\textbf{I}$ is the unit matrix, $\zeta$ is a small number,
and takes $10^{-14}$ in the computation. The matrix $(\zeta
\textbf{I}+\textbf{A}^T\textbf{A})$ is invertible, and the quadratic
polynomial can be determined uniquely. Meanwhile, if the matrix
$\textbf{A}^T\textbf{A}$ is invertible, the difference between the
solution of Eq.\eqref{quadratic-2} and Eq.\eqref{quadratic-1} is
$O(\zeta)$. With the coefficient $\textbf{a}=(a_1,...,a_5)^T$, the
point value of $P^2(x,y)$ at the Gaussian quadrature point
$(x_G,y_G)$ can be written as a linear combination of the cell
averaged values
\begin{equation}\label{polynomial-eta}
P^2(x_G,y_G)=\sum_{j=0}^{12}\eta_jW_{i_j},
\end{equation}
where $\eta_j=0$ if $E_j=0$.

\begin{rmk} Denote $\widetilde{A}\equiv \textbf{A}^T\textbf{A}$, and the
inverse of $(\zeta I+\widetilde{A})$ can be expressed as
\begin{equation}\label{proof3}
(\zeta I+\widetilde{A})^{-1}=c_0I+c_1\widetilde{A}+c_2\widetilde{A}^2+...+c_{k-1}\widetilde{A}^{k-1}.
\end{equation}
Multiplying $\widetilde{A}$ on both side of Eq.\eqref{proof3},  we have
\begin{equation*}
c_0\zeta I+(c_0+c_1\zeta)\widetilde{A}+(c_1+c_2\zeta)\widetilde{A}^2+...+(c_{k-2}+c_{k-1}\zeta)\widetilde{A}^{k-1}+c_{k-1}\widetilde{A}^{k}=I,
\end{equation*}
According to the Hamilton-Caylay theorem, $\widetilde{A}^{k}$ can be
expressed as
\begin{equation}\label{proof4}
\widetilde{A}^{k}=q_{0}I +q_{1}\widetilde{A}+...+q_{k-2}\widetilde{A}^{k-2} +q_{k-1}\widetilde{A}^{k-1},
\end{equation}
where $q_{n},~n=0,1,2...,k-1$ are given by
\begin{align*}
\left| \lambda I-\widetilde{A} \right|=-q_{0} -q_{1}\lambda -... -q_{k-2}\lambda^{k-2} -q_{k-1}\lambda^{k-1} +\lambda^k.
\end{align*}
Substituting Eq.\eqref{proof4} into Eq.\eqref{proof3}, the linear
system for $(c_0,c_1,...,c_{k-1})$ can be written as
\begin{equation*}
\left(
  \begin{array}{ccccc}
    \zeta &       &         &        & q_0             \\
    1     & \zeta &         &        & q_1             \\
          & 1     & \zeta   &        & q_2             \\
          &       & \ddots  & \ddots & \vdots          \\
          &       &         & 1      & (\zeta+q_{k-1}) \\
  \end{array}
\right) \left(
  \begin{array}{c}
    c_0 \\
    c_1 \\
    c_2 \\
    \vdots \\
    c_{k-1}
  \end{array}
\right)= \left(
  \begin{array}{c}
    1 \\
    0 \\
    0 \\
    \vdots \\
    0
  \end{array}
\right).
\end{equation*}
If the matrix $\textbf{A}^T\textbf{A}$ is invertible, $q_0\neq 0$
and the solution of the linear system above can be given as
\begin{align*}
\textbf{c}=\textbf{c}^\ast+O(\zeta)(1,1,...,1).
\end{align*}
where $\textbf{c}^\ast$ corresponds to the solution of the system
with $\zeta=0$. Substituting $\textbf{c}$ into Eq.\eqref{proof3},
and supposing the solution of Eq.\eqref{quadratic-1} is
$\textbf{a}^\ast$, the solution of Eq.\eqref{quadratic-2} can be
given as
\begin{align*}
\textbf{a}^\ast=\textbf{a}+O(\zeta)(1,1,...,1).
\end{align*}
Thus, the difference between solutions of Eq.\eqref{quadratic-2} and
Eq.\eqref{quadratic-1} is the order of $O(\zeta)$.
\end{rmk}

Similar with the standard WENO reconstruction
\cite{WENO-Liu,Hu-Shu}, twelve sub-stencils $S_{j}, j=1,...,12$ are
selected from the large stencil given in Figure.\ref{quadstencil}
\begin{align*}
P_{1}^1 ~&\text{on}~ S_1=\{i_0,i_1,i_2\}, ~~~P_{2}^1 ~\text{on}~ S_2=\{i_0,i_2,i_3\}, ~~~P_{3}^1 ~\text{on}~ S_3=\{i_0,i_3,i_4\},\\
P_{4}^1 ~&\text{on}~ S_4=\{i_0,i_4,i_1\}, ~~~P_{5}^1 ~\text{on}~ S_5=\{i_0,i_1,i_5\}, ~~~P_{6}^1 ~\text{on}~ S_6=\{i_0,i_1,i_6\},\\
P_{7}^1 ~&\text{on}~ S_7=\{i_0,i_3,i_7\}, ~~~P_{8}^1 ~\text{on}~ S_8=\{i_0,i_3,i_8\}, ~~~P_{9}^1 ~\text{on}~ S_9=\{i_0,i_4,i_9\},\\
P_{10}^1 ~&\text{on}~ S_{10}=\{i_0,i_2,i_{10}\}, P_{11}^1
~\text{on}~ S_{11}=\{i_0,i_2,i_{11}\}, P_{12}^1 ~\text{on}~
S_{12}=\{i_0,i_4,i_{12}\}.
\end{align*}
Twelve candidate linear polynomials $P_{j}^1$ corresponding to the
candidate sub-stencils can be constructed as follows
\begin{align*}
\frac{1}{\left|
\Omega_{S_{jk}}\right|}\int_{\Omega_{S_{jk}}}P^1_j(x,y)\text{d}x\text{d}y=W_{S_{jk}}, ~k=0,1,2.
\end{align*}
where $\Omega_{S_{jk}}$ is the $(k+1)$-th cell in the sub-stencils
$S_{j}$, and $W_{S_{jk}}$ is the cell averaged value over cell
$\Omega_{S_{jk}}$. With the robust least square method introduced
above, twelve linear polynomials can be fully determined. Similarly,
the point value of $P_{j}^1$ at the Gaussian quadrature point
$(x_G,y_G)$ can be written as a linear combination of the cell
averaged values $W_{jk}$. However, in order to calculate the linear
weights, it must be taken into account that the cell $i_5,...i_{12}$
may not exist for arbitrary unstructured quadrilateral meshes.
According to the local topology, the candidate polynomials is
rewritten as follows
\begin{enumerate}
\item For $j=1,...,8$, the linear polynomial can be expressed as
\begin{equation*}
P^1_j(x_G,y_G)=\sum_{k=0}^{2}b_{jk}W_{S_{jk}}.
\end{equation*}
\item For $j=9,...,12$, the linear polynomial can be expressed as
\begin{equation*}
P^1_j(x_G,y_G)=\sum_{k=0}^{1}b_{jk}W_{S_{jk}}+E_j \cdot
b_{j,2}W_{S_{j,2}}+(1-E_j) \cdot b_{j,2}W_{S_{j-4,2}}.
\end{equation*}
\end{enumerate}
For the linear scheme, the linear combination of $P_{j}^1$
\begin{align*}
R(x,y)=\sum_{j=1}^{12}\gamma_j P_{j}^1(x,y),
\end{align*}
satisfies
\begin{equation}\label{linRelation}
R(x_G,y_G)=P^2(x_G,y_G),
\end{equation}
where $(x_G, y_G)$ is the Gaussian quadrature point, and $\gamma_j$
is the linear weights only depending on the local geometry of mesh.
Substituting $P^1_j(x,y)$ into $R(x,y)$ at $(x_G,y_G)$, we have
\begin{equation} \label{lincombination}
R(x_G,y_G)=\sum_{j=0}^{12}\sum_{k=1}^{12}b_{jk}\gamma_k W_{i_j},
\end{equation}
Comparing the coefficient of the cell average value $W_{i_j}$ in
Eq.\eqref{polynomial-eta} and Eq.\eqref{lincombination}, the
following thirteen linear equations can be obtained
\begin{equation*}
\sum_{j=1}^{12}b_{kj}\cdot\gamma_j=\eta_k,~k=0,...,12,
\end{equation*}
and its matrix form is
\begin{equation} \label{linearweight-matrix}
\textbf{B}\cdot\boldsymbol{\gamma}=\boldsymbol{\eta},
\end{equation}
where $\textbf{B}=(b_{j,k})_{13\times12}$,
$\boldsymbol{\gamma}=(\gamma_1,\gamma_2,...,\gamma_{12})^T$, and
$\boldsymbol{\eta}=(\eta_0,\eta_1,...,\eta_{12})^T$. Similar with
the analysis in \cite{Hu-Shu}, the linear system for
Eq.\eqref{linearweight-matrix} contains thirteen linear equations,
and generally the system is under-determined with $\text{Rank}=11$.
$P_{j}^1(x,y)$ and $P^2(x,y)$ reproduce the linear function exactly
and Eq.\eqref{linRelation} is valid for $W(x,y)=1$, $W(x,y)=x$ and
$W(x,y)=y$ under the following identical constraint
\begin{equation*}
\displaystyle\sum_{j=1}^{12}\gamma_j=1.
\end{equation*}
Two equations can be eliminated from the linear system with thirteen
equations. The traditional least square becomes under-determined,
and the following least square method are used to obtain the linear
weights
\begin{equation*}
(\zeta \textbf{I}+\textbf{B}^T\textbf{B})\cdot
\boldsymbol{\gamma}=\textbf{B}^T\boldsymbol{\eta}.
\end{equation*}

\subsection{Linear reconstruction on triangular meshes}
\begin{figure}[!htb]
\centering
\includegraphics[width=0.45\textwidth]{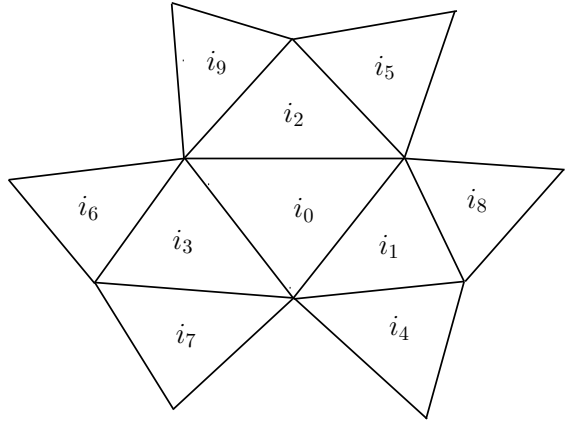}
\caption{\label{tri-stencil} Stencil for third-order scheme on
triangular meshes.}
\end{figure}

For the reconstruction on unstructured triangular meshes, the
procedure is similar with the reconstruction on quadrilateral
meshes, except the choice for stencil and sub-stencils. The stencil
is given in Figure.\ref{tri-stencil}, and the quadratic polynomial
Eq.\eqref{qua-def} can be constructed on cell $\Omega_{i_0}$ as
well. For the cell $i_0$, the neighboring cell $i_1, i_2, i_3$ must
exist and their neighboring cells $i_4, i_5, ..., i_9$ may not
exist, and the similar preconditioning is needed as well, and the
procedure is given as follows
\begin{enumerate}
\item initialize the parameter $E_j=1$, which is corresponding to the cell $i_j,j=1,...,9$;
\item read the neighboring cell $i_1, i_2, i_3$ of cell $i_0$ from the mesh file;
\item read the neighboring cell $i_j$ of cell $i_{j-3},~j=4,5,6$, if each
of them is not exist, record $E_j=0$ and $i_{j}$ is replaced by
$i_0$;
\item read the neighboring cell $i_j$ of cell $i_{j-3},~j=7,8,9$, if each
of them is not exist, record $E_j=0$ and $i_{j}$ is replaced by
$i_{j-3}$.
\end{enumerate}

In order to obtain the quadratic polynomial, the following
constrains need to be satisfied
\begin{equation*}
\frac{1}{\left| \Omega_{i_j}
\right|}\int_{\Omega_{i_j}}P^2(x,y)\text{d}x\text{d}y=W_{i_j},~ j=1,...,9.
\end{equation*}
and the matrix form can be written as
\begin{equation*}
\textbf{A}\cdot \textbf{a}=\Delta\textbf{W},
\end{equation*}
where $\textbf{a}=(a_1,...,a_5)^T$, $\Delta\textbf{W}=(W_{i_1}-W_{i_0},...,W_{i_{9}}-W_{i_0})^T$,
and $\textbf{A}=(A_{jk}),~j=1,...,9,~k=1,...,5$ is given as
\begin{equation*}
A_{jk}=E_j\cdot \frac{1}{\left| \Omega_{i_j}
\right|}\displaystyle\iint_{\Omega_{i_j}}p^k(x,y)\text{d}x\text{d}y.
\end{equation*}
Similar with the quadrilateral meshes, the polynomial $P^2(x,y)$ can
be determined. For the triangular meshes, nine candidate
sub-stencils $S_{j}, j=1,...,9$ are selected
\begin{align*}
&P_1^1 ~\text{on}~ S_1=\{i_0,i_1,i_2\}, ~ P_2^1 ~\text{on}~ S_2=\{i_0,i_2,i_3\}, ~ P_3^1 ~\text{on}~ S_3=\{i_0,i_3,i_1\},\\
&P_4^1 ~\text{on}~ S_4=\{i_0,i_1,i_4\}, ~ P_5^1 ~\text{on}~ S_5=\{i_0,i_2,i_5\}, ~ P_6^1 ~\text{on}~ S_6=\{i_0,i_3,i_6\},\\
&P_7^1 ~\text{on}~ S_7=\{i_0,i_3,i_7\}, ~ P_8^1 ~\text{on}~ S_8=\{i_0,i_1,i_8\}, ~ P_9^1 ~\text{on}~ S_9=\{i_0,i_2,i_9\},
\end{align*}
and nine linear polynomials can be fully determined. Similarly, the
point value of $P_{j}^1$ at the Gaussian quadrature point
$(x_G,y_G)$ can be written as a linear combination of the cell
averaged values $W_{jk}$. However, the cell $i_4,...i_9$ may not
exist for arbitrary unstructured meshes. In order to calculate
linear weights later, the candidate polynomials is given according
to the local topology as follows
\begin{enumerate}
\item For $j=1,...,6$, the linear polynomial can be expressed as
\begin{equation*}
P^1_j(x_G,y_G)=\sum_{k=0}^{2}b_{jk}W_{S_{jk}}.
\end{equation*}
\item For $j=7,...,9$, the linear polynomial can be expressed as
\begin{equation*}
P^1_j(x_G,y_G)=\sum_{k=0}^{1}b_{jk}W_{S_{jk}}+E_j \cdot
b_{j,2}W_{S_{j,2}}+(1-E_j) \cdot b_{j,2}W_{S_{j-3,2}},
\end{equation*}
\end{enumerate}
For the linear scheme, the linear combination
\begin{align*}
R(x,y)=\sum_{j=1}^{9}\gamma_j P_{j}^{1}(x,y),
\end{align*}
satisfies
\begin{equation}\label{combination}
R(x_G,y_G)=P^2(x_G,y_G),
\end{equation}
where $(x_G, y_G)$ is the Gaussian quadrature point, and $\gamma_j$
is the linear weights. The procedure of determining the linear
weights for quadrilateral meshes can be extended to triangular
meshes directly.

\subsection{Optimization approach for linear weights}
\begin{figure}[!htb]
\centering
\includegraphics[width=0.475\textwidth]{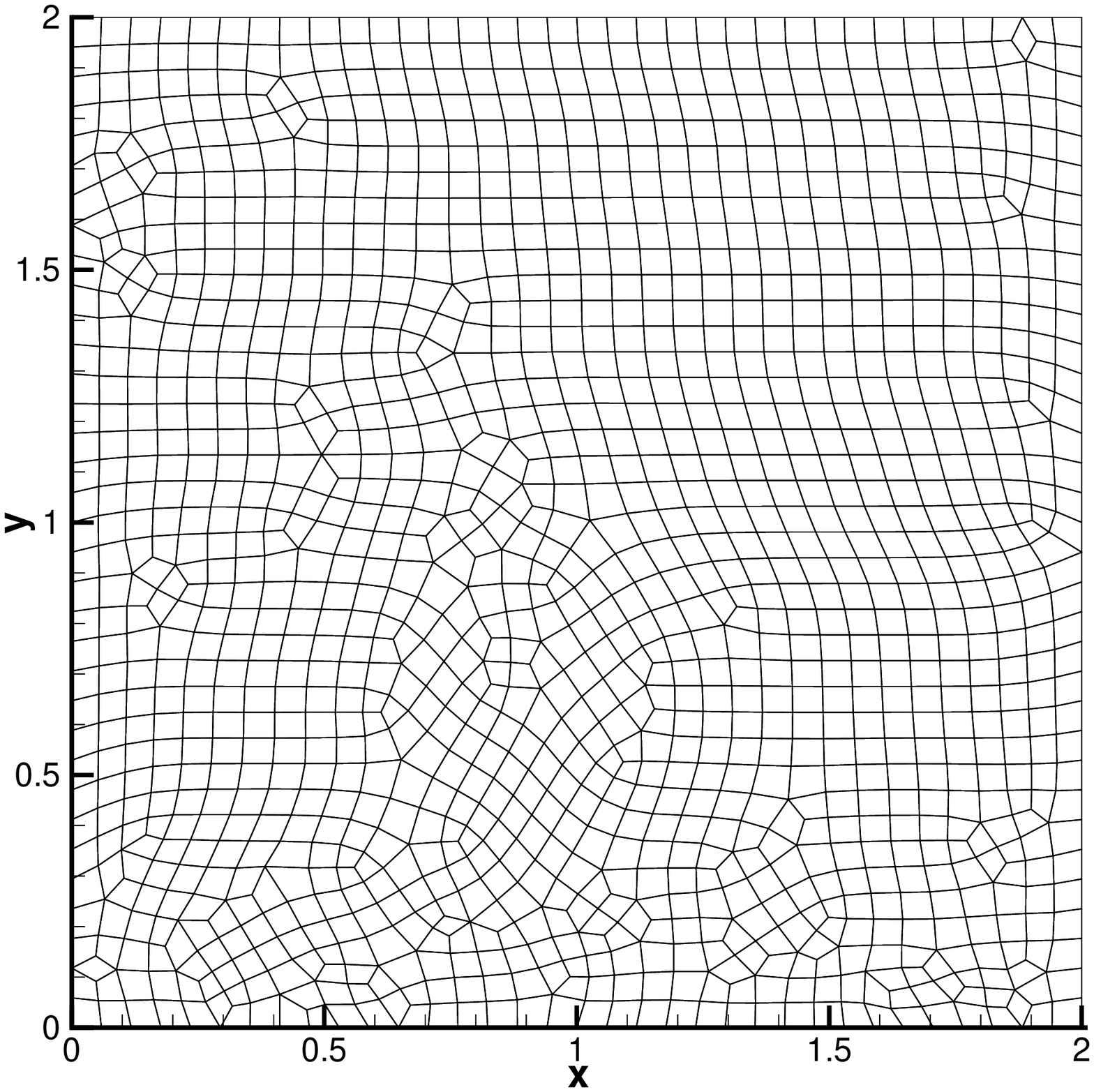}
\includegraphics[width=0.475\textwidth]{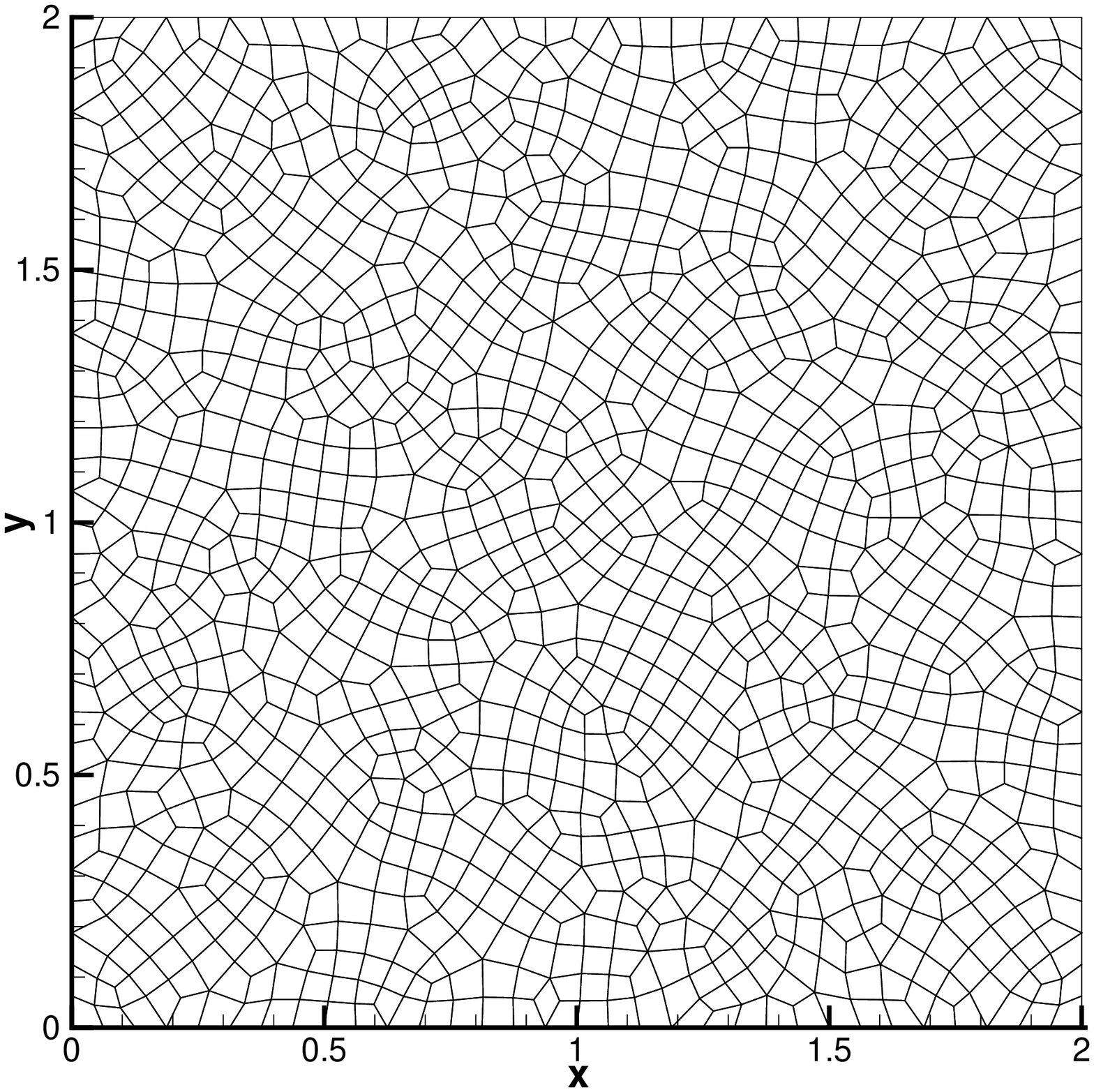}
\caption{\label{Mesh-Accuracy-1} The quadrilateral meshes with
cell size $h=1/16$: regular mesh (left) and irregular mesh (right).}
\centering
\includegraphics[width=0.475\textwidth]{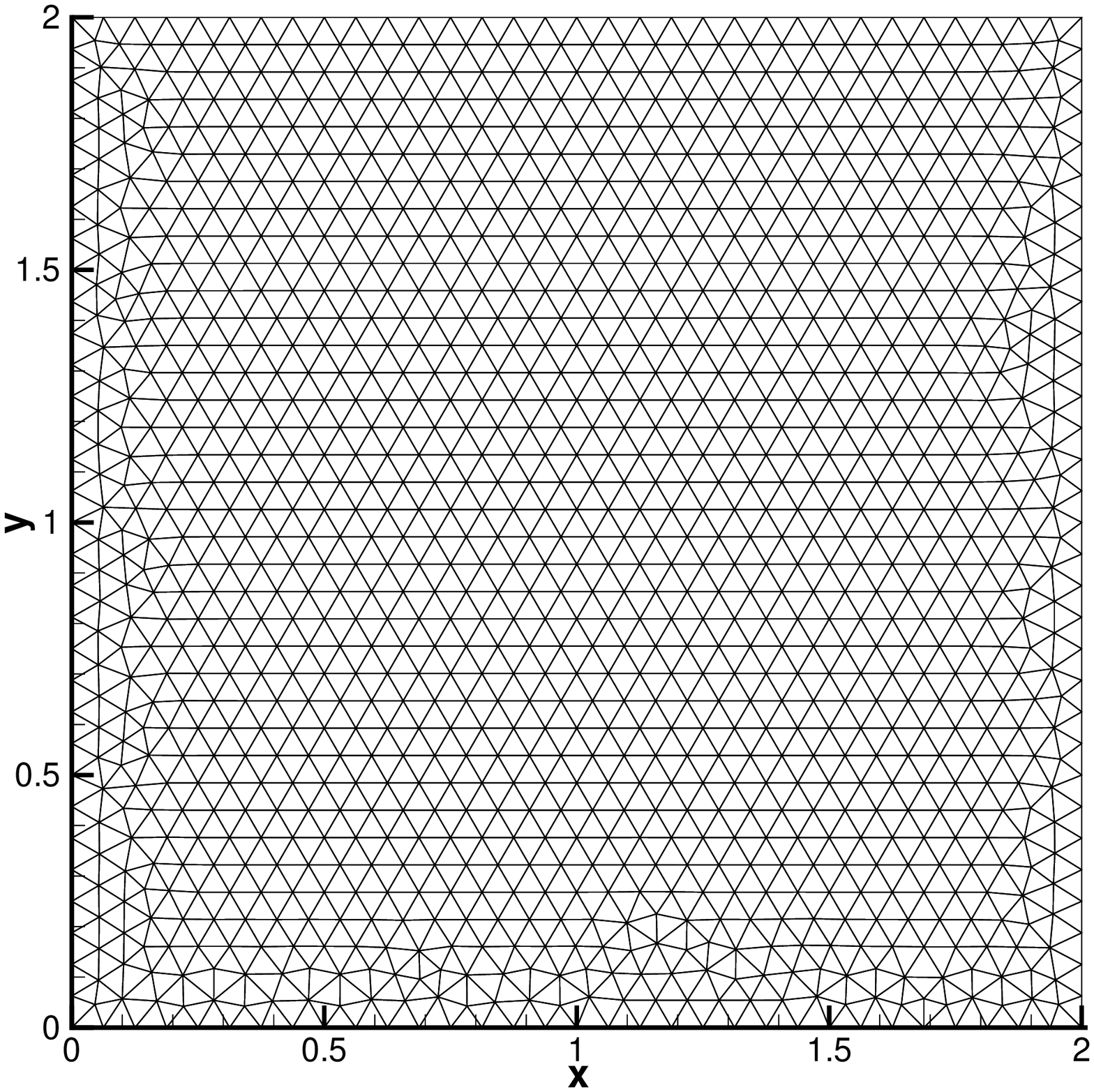}
\includegraphics[width=0.475\textwidth]{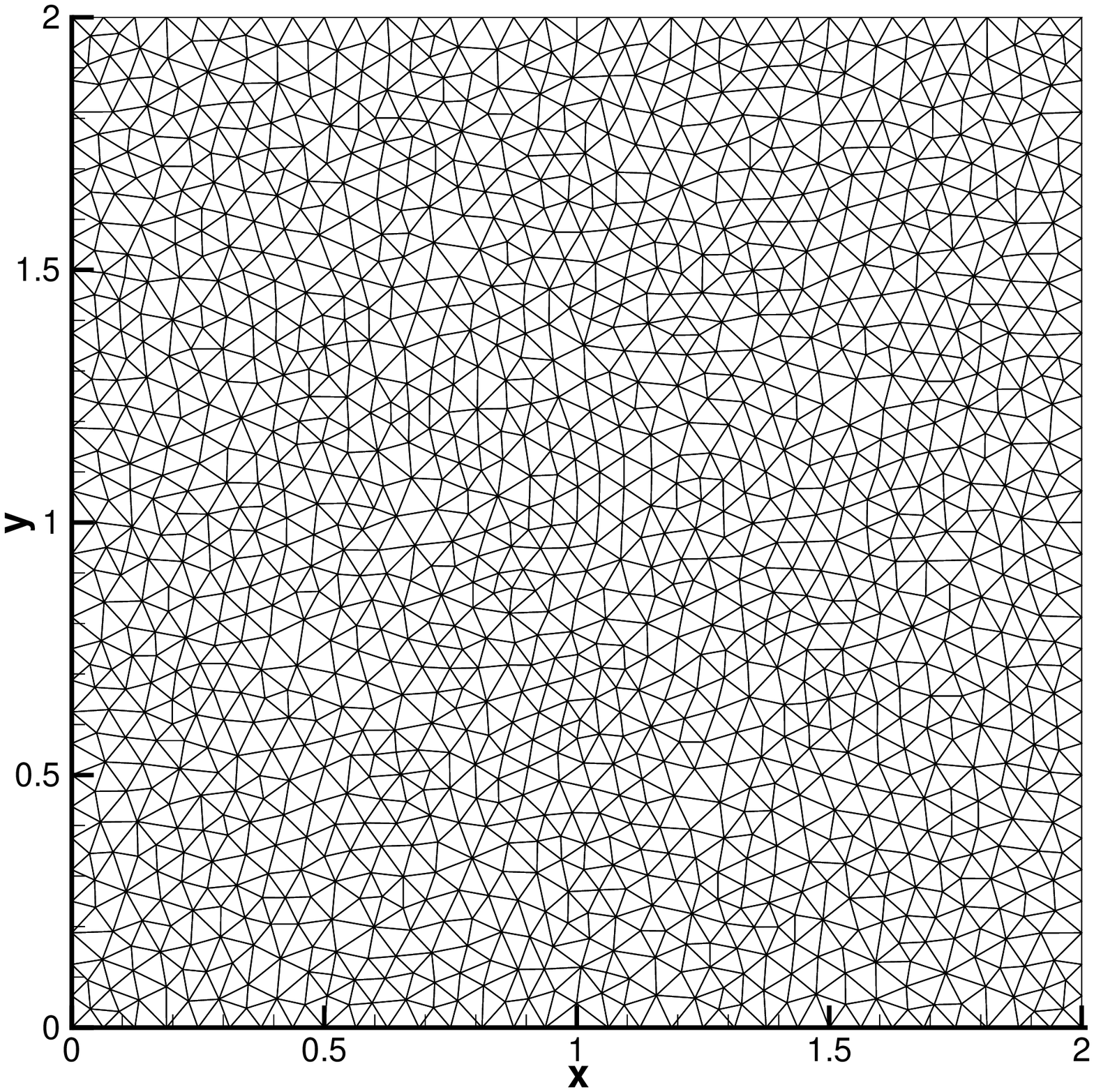}
\caption{\label{Mesh-Accuracy-2} The triangular meshes for with
cell size $h=1/16$: regular mesh (left) and irregular mesh (right).}
\end{figure}

With the procedure introduced above, the linear weights for WENO
reconstruction can be obtained. To deal with the flow with
discontinuities, the non-linear weights need to be constructed,
which is one of the most important step in the WENO schemes.
However, for the unstructured meshes, the very large linear weights
and negative linear weights appear commonly, which makes the WENO
schemes unstable even for the smooth flow problems. In order to
construct a robust WENO schemes on unstructured meshes, the very
large linear weights and negative linear weights need to be treated
carefully. To illustrate the very large weights, the quadrilateral
and triangular meshes in Figure.\ref{Mesh-Accuracy-1} and
Figure.\ref{Mesh-Accuracy-2} are considered, where the meshes with
cell size $h=1/16$ are given. The meshes in the left column are
closer to the uniform quadrilateral and triangular meshes, and they
denoted as the regular meshes. Meanwhile, the meshes in the right
column are denoted as irregular meshes. The maximum of the liner
weights $\left|\gamma\right|_{max}$ for these meshes are given in
Table.\ref{tab-mesh-1}, where the maximum becomes larger with the
mesh refinement. Such very large weights will trigger the
instability of WENO scheme.

\begin{table}[!h]
\begin{center}
\def\temptablewidth{0.8\textwidth}
{\rule{\temptablewidth}{1.0pt}}
\begin{tabular*}{\temptablewidth}{@{\extracolsep{\fill}}c|c|c|c|c}
mesh  size & regular quad & irregular quad &  regular tri & irregular tri   \\
\hline
$1/8$  & 4.3589E01       &  4.9770E02  & 2.1023E02       & 2.6012E02  \\
$1/16$ & 2.9415E02       &  8.4584E02  & 1.0262E02       & 2.9049E05  \\
$1/32$ & 1.9161E02       &  2.7036E03  & 9.6171E02       & 2.6894E03  \\
$1/64$ & 6.6936E02       &  2.7012E04  & 1.9599E03       & 1.1906E04
\end{tabular*}
{\rule{\temptablewidth}{1.0pt}}
\caption{\label{tab-mesh-1}Maximum
of the liner weights $\left|\gamma\right|_{max}$  for the
quadrilateral and triangular meshes.}
\end{center}
\begin{center}
\def\temptablewidth{0.8\textwidth}
{\rule{\temptablewidth}{1.0pt}}
\begin{tabular*}{\temptablewidth}{@{\extracolsep{\fill}}c|c|c|c|c}
mesh  size & regular quad & irregular quad &  regular tri & irregular tri   \\
\hline
$1/8$  &  1.3897E00   & 1.4613E00   & 1.6410E00  & 2.0553E00 \\
$1/16$ &  1.3719E00   & 1.5690E00   & 1.7563E00  & 2.0730E00 \\
$1/32$ &  1.4029E00   & 1.5359E00   & 1.6214E00  & 8.5631E00 \\
$1/64$ &  1.3910E00   & 1.9566E00   & 1.6269E00  & 6.4324E00
\end{tabular*}
{\rule{\temptablewidth}{1.0pt}}
\caption{\label{tab-mesh-2} Maximum
of the optimal liner weights $\left|\gamma\right|_{max}^{opt}$  for
the quadrilateral and triangular meshes.}
\end{center}
\end{table}

The meshes with lower quality will result in the very large linear
weights, and an optimization approach is introduced to deal with
them. In the current optimization approach, the weighting parameters
is introduced for each cell, such that the ill cell contributes
little to the quadratic polynomial in Eq.\eqref{polynomial-eta}. For
quadrilateral meshes, the weighting parameters are defined
\begin{equation*}
\begin{cases}
d_j=1, & ~~j=1,...,4, \\
d_j=\displaystyle\frac{1}{\max(1,\left| \gamma_j\right|)}, &
~~j=5,...,12.
\end{cases}
\end{equation*}
For triangular meshes, the weighting parameters are defined
\begin{equation*} 
\begin{cases}
d_j=1, j=1,...,3, \\
d_j=\displaystyle\frac{1}{\max(1,\left| \gamma_j\right|,\left|
\gamma_k\right|)}, (j,k)=(4,8),(8,4),(5,9),(9,5),(6,7),(7,6),
\end{cases}
\end{equation*}
The optimized coefficients $\textbf{a}=(a_1,...,a_5)$ of the
quadratic polynomial in Eq.\eqref{polynomial-eta} are given by the
following weighted linear system
\begin{equation*}
\sum_{k=1}^5 d_j\cdot A_{jk}a_k=d_j\cdot (W_{i_j}-W_i),~j=1,...,M,
\end{equation*}
where $M=12$ for the quadrilateral meshes, and $M=9$ for the
triangular meshes. The maximum of optimized liner weights
$\left|\gamma\right|_{max}^{opt}$ for the meshes in
Figure.\ref{Mesh-Accuracy-1} and Figure.\ref{Mesh-Accuracy-2} are
given in Table.\ref{tab-mesh-2}. With the mesh refinement, the
maximum of liner weights $\left|\gamma\right|_{max}^{opt}$ is
reduced greatly. This linear system is solved with the least square
method proposed above. With the procedure above, the maximum
$\gamma_j$ becomes the order $O(1)$.

\subsection{Nonlinear weights}
Before implementing the WENO reconstruction, the splitting technique
\cite{splitting-weights} is considered for the negative linear
weights which is obtained by the optimized approach as follows
\begin{equation*}
\widetilde{\gamma}_j^{+}=\frac{1}{2}(\gamma_j + \theta
\left|\gamma_j\right|),
~~\widetilde{\gamma}_j^{-}=\widetilde{\gamma}_j^{+}-\gamma_j,
\end{equation*}
where $\theta=3$ is taken in numerical tests, and the scaled
non-negative weights are given by
\begin{equation*}
\gamma_j^{\pm}=\frac{\widetilde{\gamma}_j^{\pm}}{\sigma^{\pm}},~\sigma^{\pm}=\sum_{j}\widetilde{\gamma}_j^{\pm}.
\end{equation*}
With the non-negative weights, the splitting polynomial is written as
\begin{equation*}
R^{\pm}(x_G,y_G)=\sum_j\gamma_j^{\pm}P^1_j(x_G,y_G),
\end{equation*}
which satisfies
\begin{equation*}
R(x_G,y_G)=\sigma^{+}R^{+}(x_G,y_G)-\sigma^{-}R^{-}(x_G,y_G).
\end{equation*}
Based on the non-negative weights obtained by splitting approach,
the nonlinear weights are obtained as
\begin{equation*}
\delta_{j}^{\pm}=\frac{\alpha_{j}^{\pm}}{\sum_{l=1}\alpha_{l}^{\pm}},~~\alpha_{j}^{\pm}=\frac{\gamma_{j}^{\pm}}{(\widetilde{\beta}_j^{\pm}+\epsilon)^{2}},
\end{equation*}
where $\epsilon$ is a small positive number,
$\widetilde{\beta_j}^{\pm}$ is a new smooth indicator defined on
unstructured meshes
\begin{equation*}
\widetilde{\beta_j}^{\pm}=\beta_j\big(1+\gamma_j^{\pm}\beta_j+(\gamma_j^{\pm} \beta_j)^2 \big),
\end{equation*}
and $\beta_{j}$ is defined as
\begin{equation*}
\beta_j=\sum_{|\alpha|=1}^{K}|\Omega|^{|\alpha|-1}\int\int_{\Omega}\big(D^{\alpha}P(x,y)\big)^2dxdy.
\end{equation*}
where $\alpha$ is a multi-index and $D$ is the derivative operator.
In the new smooth indicator, the accuracy keeps the original order
in smooth regions with $IS=h^2(1+O(h))$. In addition, the linear
weights are considered in constructing the new smooth indicator.
Suppose that $\gamma_1$ is a very large linear weight and the
corresponding $P_1^1(x,y)$ contains the discontinuity. The ratio of
nonlinear weights of $\delta_1^{\pm}$ and $\delta_j^{\pm}$ becomes
\begin{equation*}
\frac{\delta_1^{\pm}}{\delta_j^{\pm}}=\frac{\gamma_1^{\pm}}{\gamma_j^{\pm}}\cdot \frac{\big(\widetilde{\beta_j}^{\pm}\big)^2}{\big(\widetilde{\beta_1}^{\pm}\big)^2}
=\frac{\gamma_1^{\pm}}{\gamma_j^{\pm}}\cdot \frac{h^2(1+O(h))}{\big( \gamma_1^{\pm}(\gamma_1^{\pm}+O(1)) \big)^2}
=\frac{h^2}{\big( \big(\gamma_1^{\pm}\big)^{1/2}(\gamma_1^{\pm}+O(1)) \big)^2}.
\end{equation*}
The contribution of the candidate polynomial $P_1^1(x,y)$ becomes
very little, and the essentially non-oscillatory property retains
for the very large linear weight. The final reconstructed value at
the Gaussian quadrature points can be written as
\begin{equation*}
R(x_G,y_G)=\sum_{j=1}(\delta_j^{+}P^1_j(x_G,y_G)-\delta_j^{-}P^1_j(x_G,y_G)).
\end{equation*}

\section{Numerical tests}
In this section, the numerical tests will be presented to validate
the current WENO scheme. For the temporal discretization, the
classical third-order TVD Runge-Kutta method \cite{TVD-RK} is used
for Eq.\eqref{semidiscrete}, which can be written as follows
\begin{align*}
W^{(1)}_i&=W^{n}+\Delta t\mathcal{L}(W^{n}_i), \\
W^{(2)}_i&=\frac{3}{4}W^{n}+\frac{1}{4}W^{(1)}+\frac{1}{4}\Delta t\mathcal{L}(W^{(1)}_i), \\
W^{n+1}_i&=\frac{1}{3}W^{n}+\frac{2}{3}W^{(2)}+\frac{2}{3} \Delta
t\mathcal{L}(W^{(2)}),
\end{align*}
where $\mathcal{L}$ is the operator corresponding to the spatial
discretization.

In order to eliminate the spurious oscillation and improve the
stability, the reconstruction for Euler system can be performed on
the characteristic variables. The Jacobian matrix at one quadrature
point is $n_x(\partial F/\partial W)_{Q}+n_y(\partial G/\partial
W)_{Q}$, and the Roe's mean matrix is used, which can be
diagnoalized by the right eigenmatrix $R$. The variables for
reconstruction are defined as $U=R^{-1}W$. With the reconstructed
values, the conservative variables can be obtained by the inverse
projection.

\begin{table}[!h]
\begin{center}
\def\temptablewidth{0.75\textwidth}
{\rule{\temptablewidth}{1.0pt}}
\begin{tabular*}{\temptablewidth}{@{\extracolsep{\fill}}c|c|cc|cc}
scheme   & mesh   &  $L^1$ error &  Order   &  $L^\infty$ error & Order \\
\hline
WENO-3    &  $1/8$  &  2.533E-01  &     ~     &   1.058E-01  &  ~     \\
regular   &  $1/16$ &  3.889E-02  &     2.70  &   2.348E-02  &  2.17  \\
meshes      &  $1/32$ &  3.730E-03  &     3.38  &   2.399E-03  &  3.29  \\
~         &  $1/64$ &  2.127E-04  &     4.13  &   1.480E-04  &  4.02  \\
\hline
WENO-3    &  $1/8$  &  2.220E-01  &     ~     &   9.791E-02  &  ~    \\
irregular &  $1/16$ &  3.262E-02  &     2.77  &   2.122E-02  &  2.21 \\
meshes     &  $1/32$ &  2.602E-03  &     3.65  &   1.877E-03  &  3.50 \\
~         &  $1/64$ &  1.544E-04  &     4.07  &   1.122E-04  &  4.06
\end{tabular*}
{\rule{\temptablewidth}{1.0pt}}
\end{center}
\vspace{-5mm}\caption{\label{tab-accuracy1} Accuracy tests: linear wave propagation on quadrilateral with WENO scheme.}
\begin{center}
\def\temptablewidth{0.75\textwidth}
{\rule{\temptablewidth}{1.0pt}}
\begin{tabular*}{\temptablewidth}{@{\extracolsep{\fill}}c|c|cc|cc}
Schemes    & mesh   & $L^1$ error &  Order ~& $L^\infty$ error & Order  \\
\hline
WENO-3     &  $1/8$  &  6.140E-02  &  ~      &   3.372E-02  &  ~     \\
regular    &  $1/16$ &  6.901E-03  &  3.15   &   4.591E-03  &  2.88  \\
meshes       &  $1/32$ &  4.822E-04  &  3.84   &   3.410E-04  &  3.75  \\
~          &  $1/64$ &  5.633E-05  &  3.10   &   2.770E-05  &  3.62  \\
\hline
WENO-3     &  $1/8$  &  5.419E-02  &   ~     &   3.211E-02  &  ~     \\
irregular  &  $1/16$ &  6.226E-03  &   3.12  &   4.505E-03  &  2.83  \\
meshes      &  $1/32$ &  4.227E-04  &   3.88  &   3.288E-04  &  3.78  \\
~          &  $1/64$ &  4.845E-05  &   3.13  &   2.550E-05  &  3.69
\end{tabular*}
{\rule{\temptablewidth}{1.0pt}}
\end{center}
\vspace{-5mm}\caption{\label{tab-accuracy2} Accuracy tests: linear wave propagation on triangular with WENO scheme.}
\end{table}

\begin{figure}[!htb]
\centering
\includegraphics[width=0.475\textwidth]{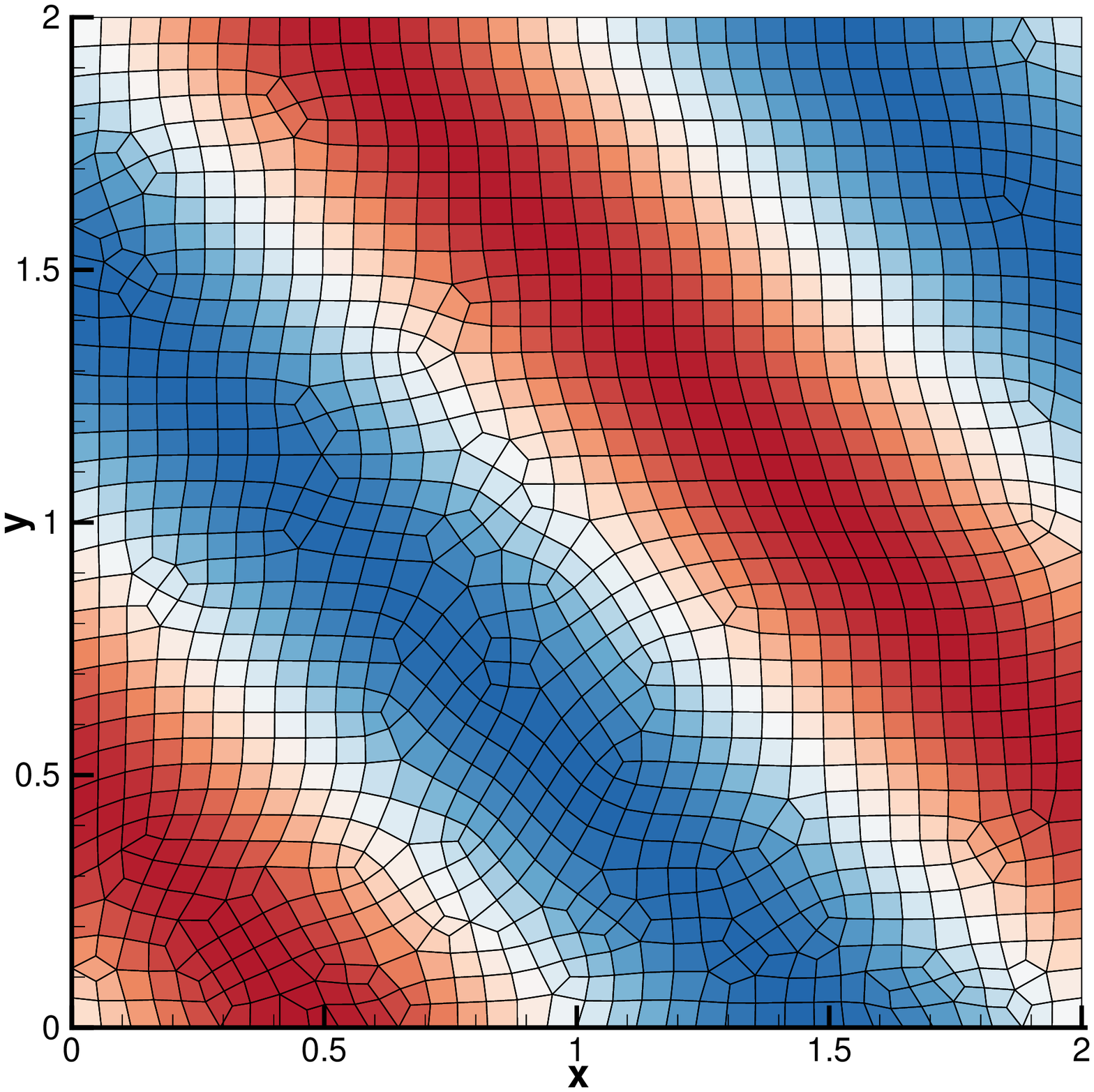}
\includegraphics[width=0.475\textwidth]{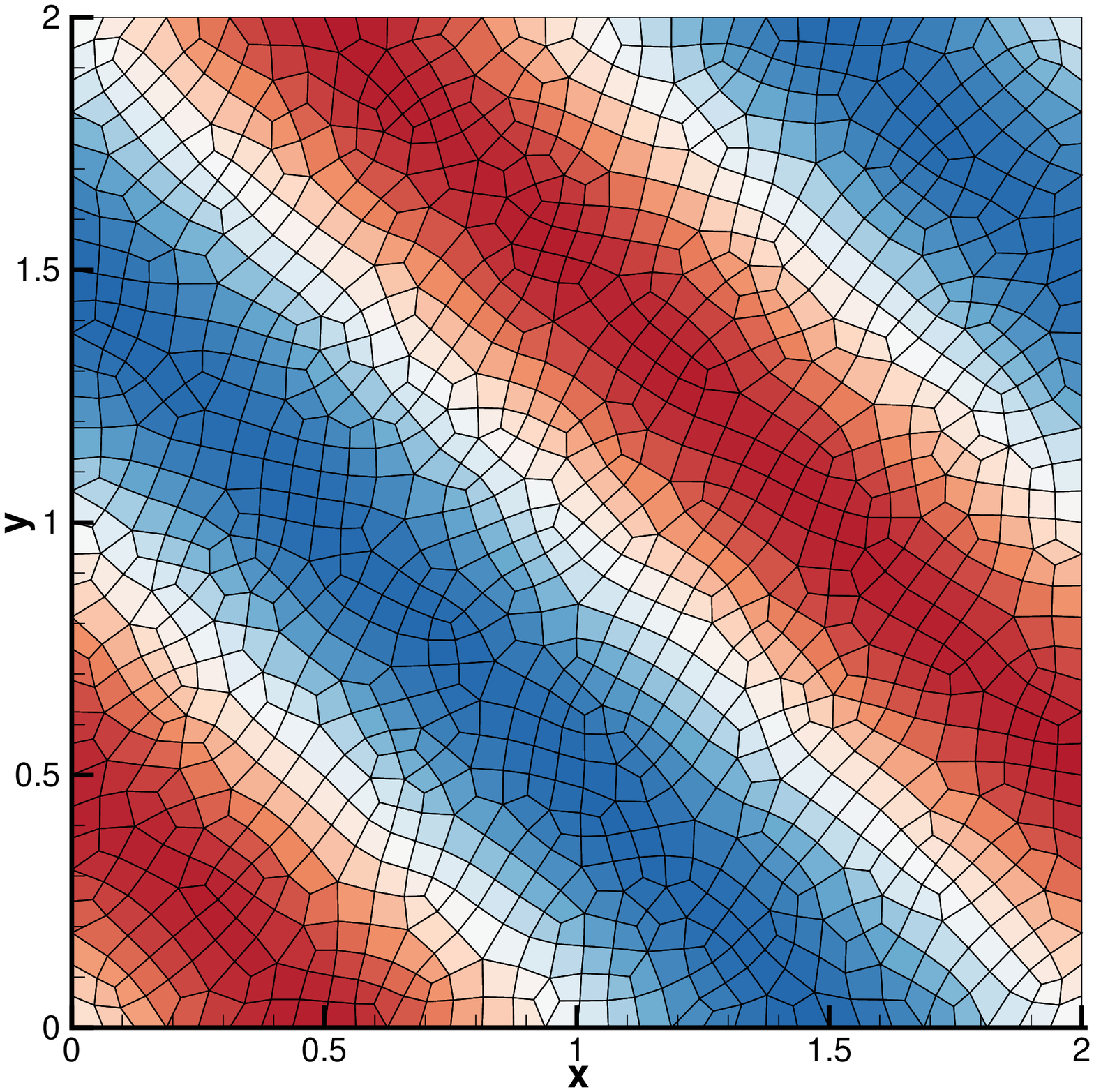}
\caption{\label{Mesh-Accuracy-1a} Advection of density perturbation: density distribution on the quadrilateral meshes with cell size $h=1/16$, regular mesh (left) and irregular mesh (right).}
\centering
\includegraphics[width=0.475\textwidth]{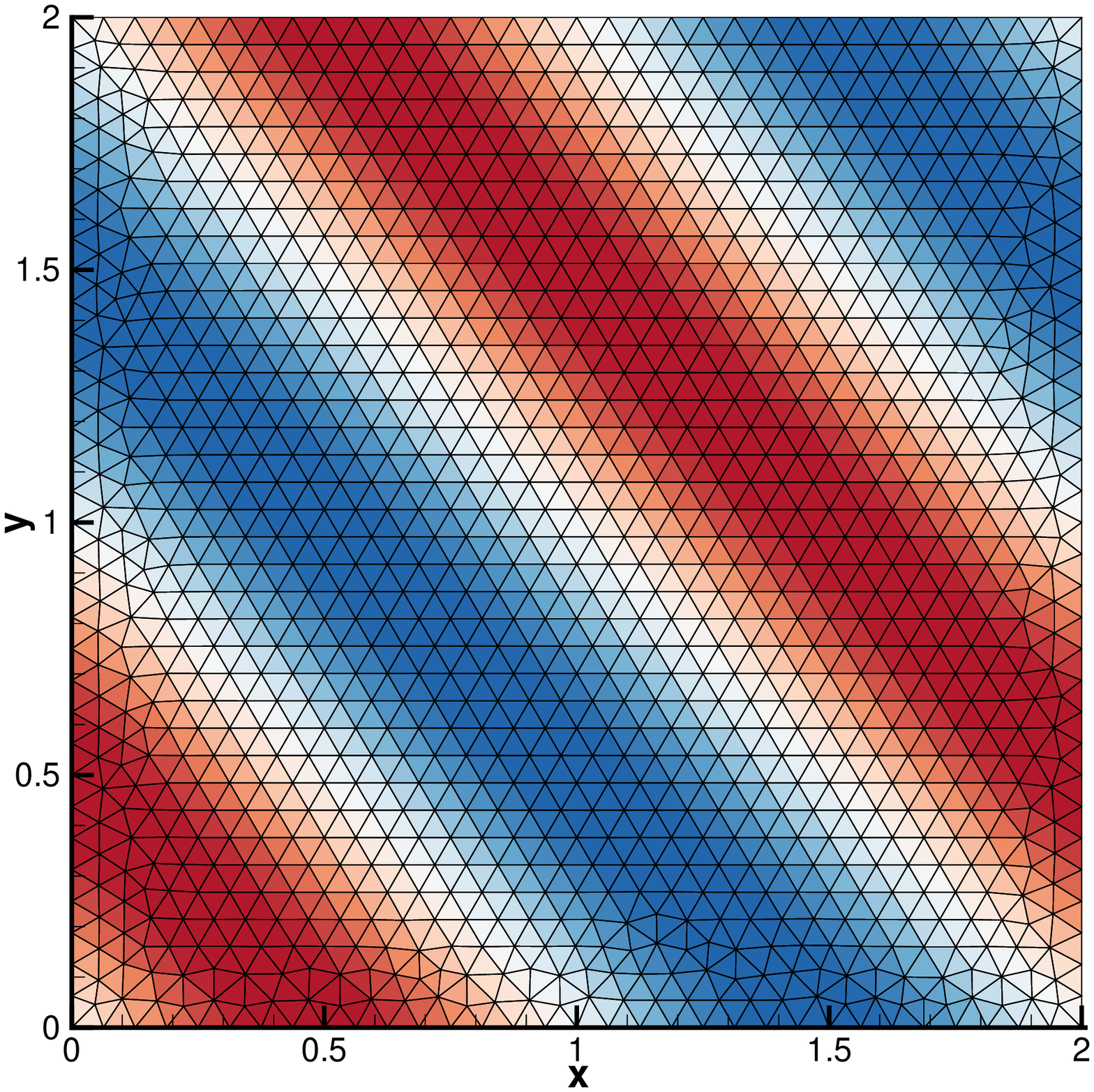}
\includegraphics[width=0.475\textwidth]{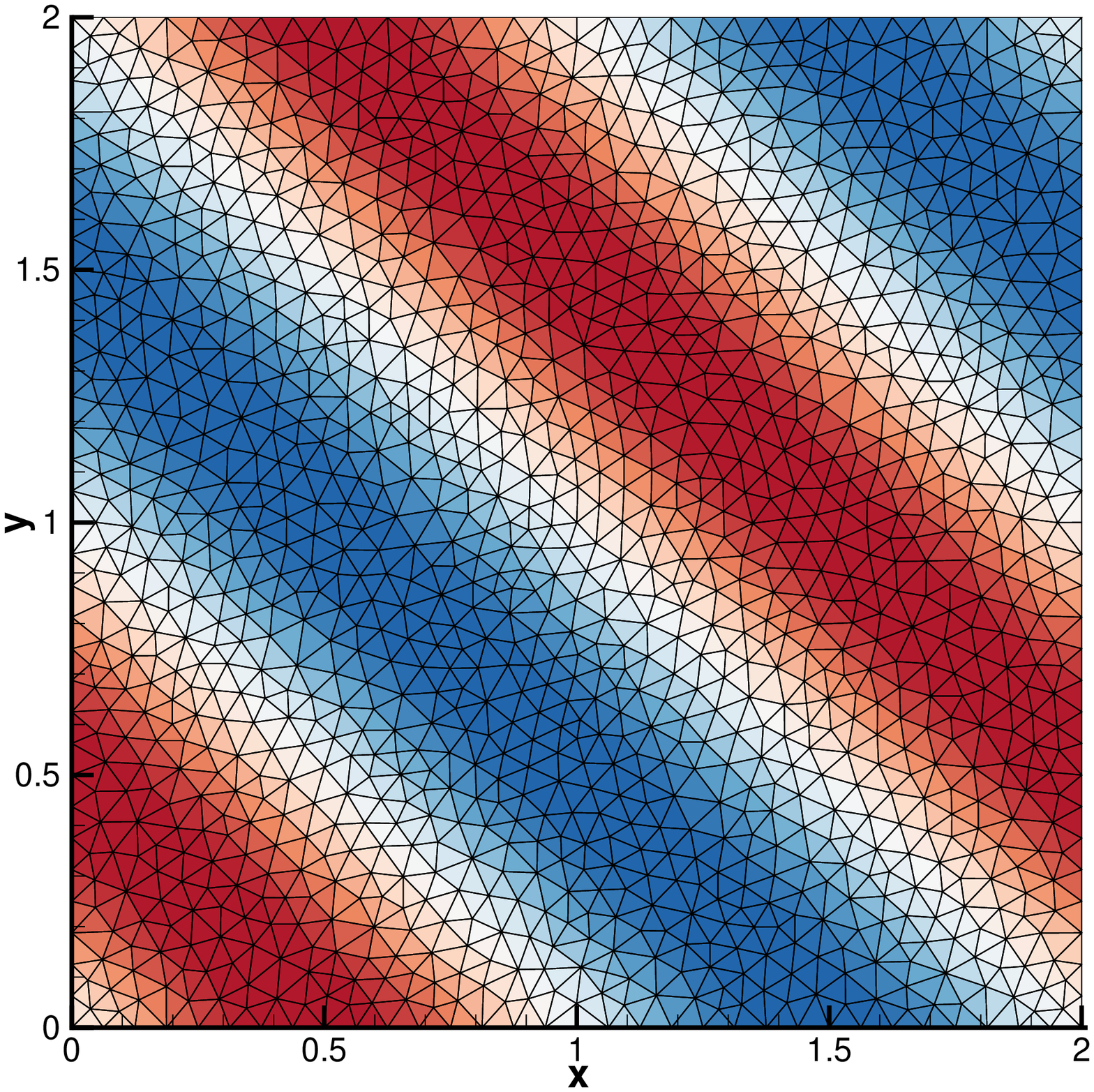}
\caption{\label{Mesh-Accuracy-2a} Advection of density perturbation: density distribution on the triangular meshes with cell size $h=1/16$, regular mesh (left) and irregular mesh (right).}
\end{figure}

\begin{table}[!h]
\begin{center}
\def\temptablewidth{0.7\textwidth}
{\rule{\temptablewidth}{1.0pt}}
\begin{tabular*}{\temptablewidth}{@{\extracolsep{\fill}}c|c|cc|cc}
scheme & mesh   &  $L^1$ error &  Order   &  $L^\infty$ error & Order \\
\hline
WENO-3    &  $1/4$  &  2.291E-01  &     ~     &   4.343E-02  &  ~     \\
regular   &  $1/8$  &  3.013E-02  &     2.93  &   6.597E-03  &  2.72  \\
meshes    &  $1/16$ &  2.984E-03  &     3.34  &   6.715E-04  &  3.30  \\
~         &  $1/32$ &  3.431E-04  &     3.12  &   6.543E-05  & 3.36\\
\hline
WENO-3    &  $1/4$  &  1.637E-01  &     ~     &   3.605E-02  &  ~    \\
irregular &  $1/8$  &  2.094E-02  &     2.97  &   5.355E-03  &  2.75 \\
meshes     &  $1/16$ &  2.346E-03  &     3.16  &   6.094E-04  &  3.14 \\
~         &  $1/32$ &  2.794E-04  &     3.07  &   6.976E-05  & 3.13
\end{tabular*}
{\rule{\temptablewidth}{1.0pt}}
\end{center}
\vspace{-5mm}\caption{\label{tab-accuracy3} Accuracy tests: isentropic vortex propagation on quadrilateral with WENO scheme.}
\begin{center}
\def\temptablewidth{0.7\textwidth}
{\rule{\temptablewidth}{1.0pt}}
\begin{tabular*}{\temptablewidth}{@{\extracolsep{\fill}}c|c|cc|cc}
Schemes    & mesh   & $L^1$ error &  Order ~& $L^\infty$ error & Order \\
\hline
WENO-3     &  $1/4$  &  5.515E-02  &  ~     &   1.236E-02  &  ~    \\
regular    &  $1/8$  &  6.256E-03  &  3.14  &   1.445E-03  &  3.10 \\
meshes     &  $1/16$ &  7.219E-04  &  3.12  &   1.368E-04  &  3.40 \\
~          &  $1/32$ &  1.328E-04  &  2.44  &   1.483E-05  &  3.21 \\
\hline
WENO-3     &  $1/4$  &  5.248E-02  &  ~     &   2.870E-02  &  ~     \\
irregular  &  $1/8$  &  5.200E-03  &  3.34  &   1.370E-03  &  4.39  \\
meshes     &  $1/16$ &  6.758E-04  &  2.94  &   1.422E-04  &  3.27  \\
~          &  $1/32$ &  1.270E-04  &  2.41  &   1.997E-05  &  2.83
\end{tabular*}
{\rule{\temptablewidth}{1.0pt}}
\end{center}
\vspace{-5mm}\caption{\label{tab-accuracy4}Accuracy tests: isentropic vortex propagation on triangular with WENO scheme.}
\end{table}

\subsection{Accuracy tests}
The first case is the advection of density
perturbation, and the initial condition is given as follows
\begin{align*}
\rho(x,y)=1+0.2\sin(\pi (x+y)),~u(x,y)=1,~v(x,y)=1,~p(x,y)=1.
\end{align*}
The computational domain is $[0,2]\times[0,2]$, and the periodic
boundary conditions are applied in both directions. To keep the
third-order accuracy, $\Delta t=0.3 h$ is used for the third-order
scheme where $h$ is the scale of meshes. In order to validate the
accuracy of the current scheme, both quadrilateral and triangular
meshes are tested.  The meshes with mesh size $h=1/16$ in Figure.\ref{Mesh-Accuracy-1}
and Figure.\ref{Mesh-Accuracy-2} are shown as example, and
the density distribution with cell size $h=1/16$ are given in Figure.\ref{Mesh-Accuracy-1a}
and Figure.\ref{Mesh-Accuracy-2a}.
The $L^1$ and $L^\infty$ errors
and convergence orders are presented in Table.\ref{tab-accuracy1} and
Table.\ref{tab-accuracy2} for quadrilateral and triangular
meshes at $t=2$.
The expected order of accuracy are obtained, and the order of error
is not affected by the quality of meshes.

The second accuracy test is the isentropic vortex propagation
problem.  The mean flow is $(\rho,u,v,p)=(1,1,1,1)$. An isentropic
vortex is added to the mean flow with the following perturbations
\begin{equation*}
(\delta u,\delta v)=\frac{\epsilon}{2\pi}e^{0.5(1-r^2)}(-\overline{y},\overline{x}),\\
\delta T=-\frac{(\gamma -1)\epsilon^2}{8\gamma\pi^2}e^{1-r^2},
\delta S=0,
\end{equation*}
where $(\overline{x},\overline{y})=(x-5,y-5)$,
$r^2=\overline{x}^2+\overline{y}^2$, and the vortex strength
$\epsilon=5$. The computational domain is $[0, 10] \times [0, 10]$, and periodic boundary condition is applied to all boundaries. The regular and irregular quadrilateral meshes, regular and irregular triangular meshes are tested respectively. The $L^1$ and $L^\infty$ errors and convergence orders for current scheme are presented in Table.\ref{tab-accuracy3} and Table.\ref{tab-accuracy4} for quadrilateral and triangular meshes at $t=2$. The expected order of accuracy are obtained, and the order of error is not affected by the quality of meshes as well. 
With the fixed mesh size $h=1/8$, the long-time evolution of the
vortex is computed. The density distribution at center line with $5$ and $10$
time periods on quadrilateral and triangular irregular meshes are given in Figure.\ref{1d-shocktube-mesh2} with mesh size $h=1/8$. The dissipation is well controlled by the current scheme even
after long time evolution. Due to more
number of cell with the same mesh scale, the scheme performs better
on triangular meshes than quadrilateral meshes.   

\begin{figure}[!htb]
\centering
\includegraphics[width=0.6\textwidth]{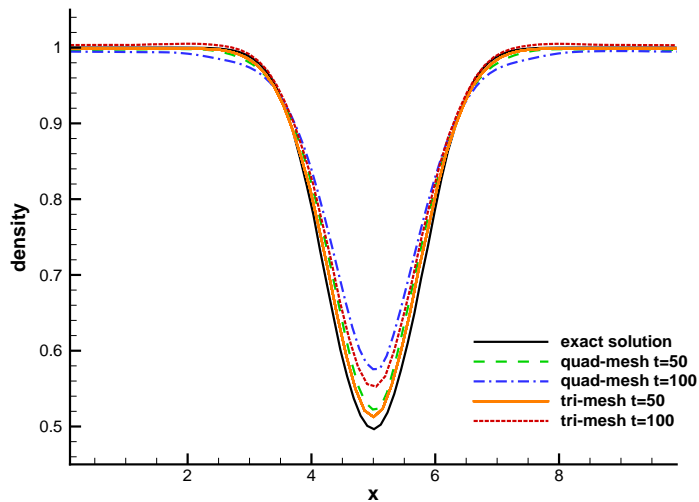}
\caption{\label{1d-shocktube-mesh2} Isentropic vortex propagation: density distribution at center line with $5$ and $10$
time periods on quadrilateral and triangular irregular meshes with mesh size $h=1/8$.}
\end{figure}

\begin{figure}[!htb]
\centering
\includegraphics[width=0.485\textwidth]{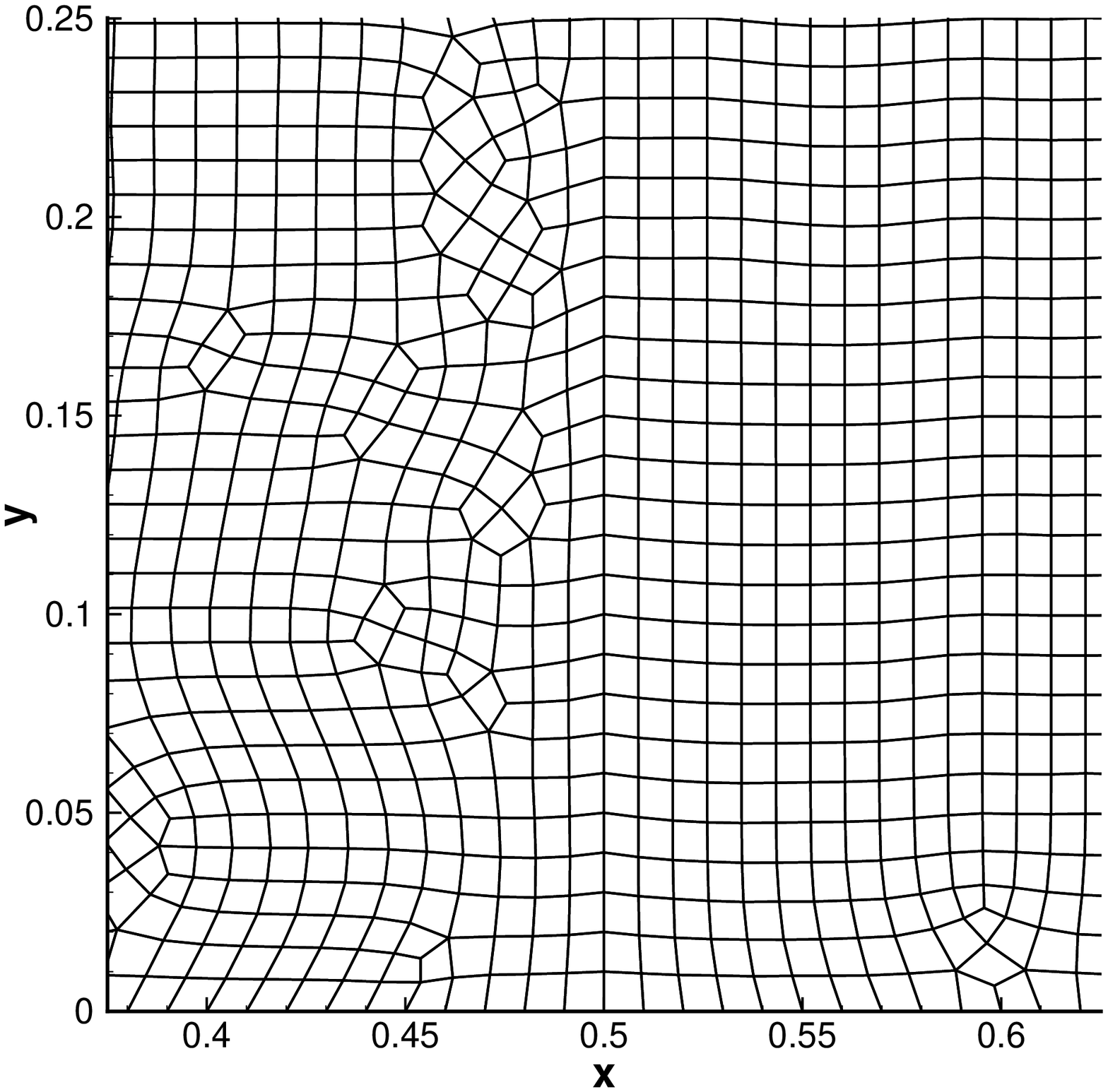}
\includegraphics[width=0.485\textwidth]{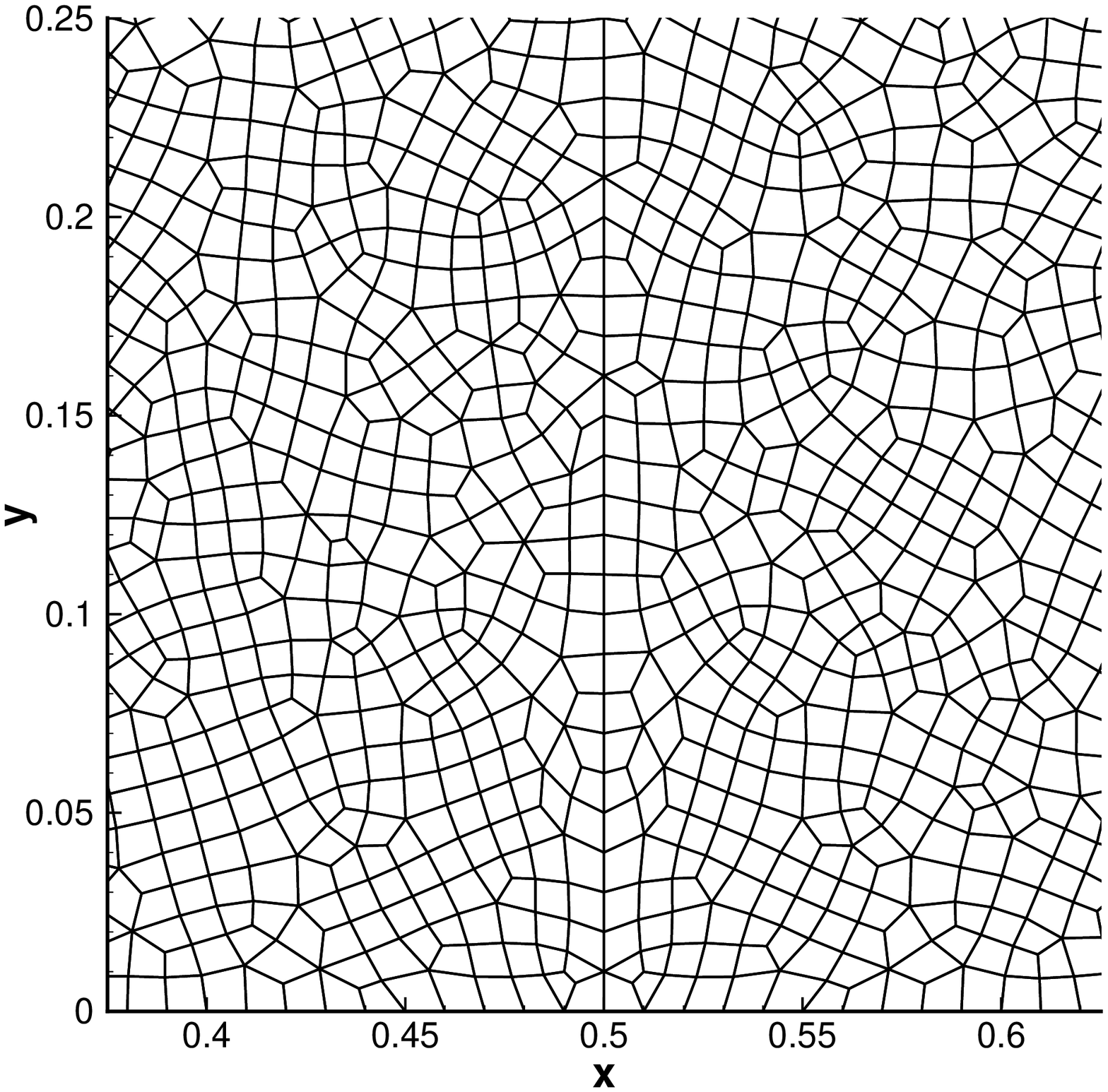}
\caption{\label{1d-shocktube-mesh1} One dimensional Riemann problem: quadrilateral regular (left) and quadrilateral irregular (right) meshes.}
\centering
\includegraphics[width=0.485\textwidth]{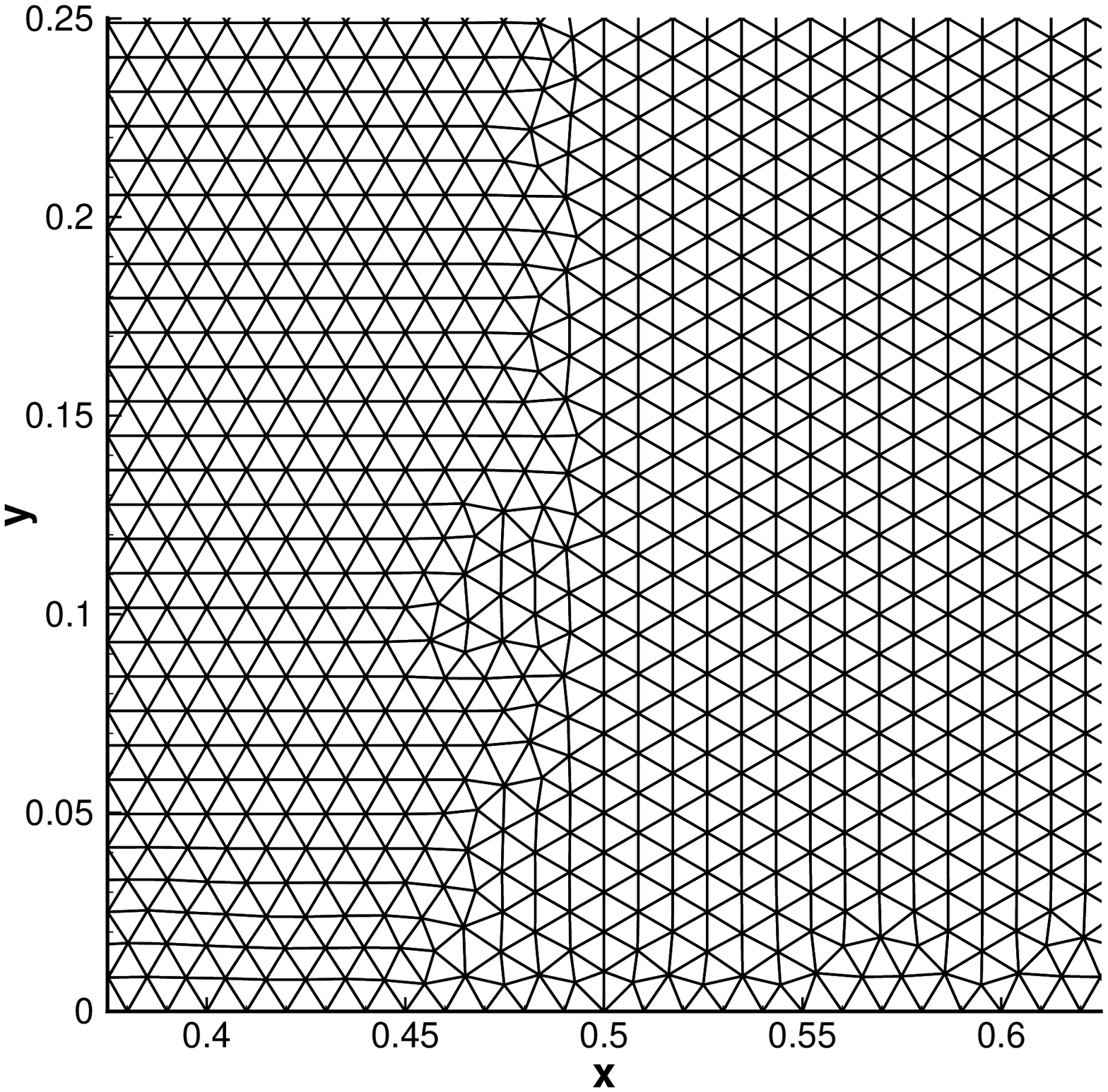}
\includegraphics[width=0.485\textwidth]{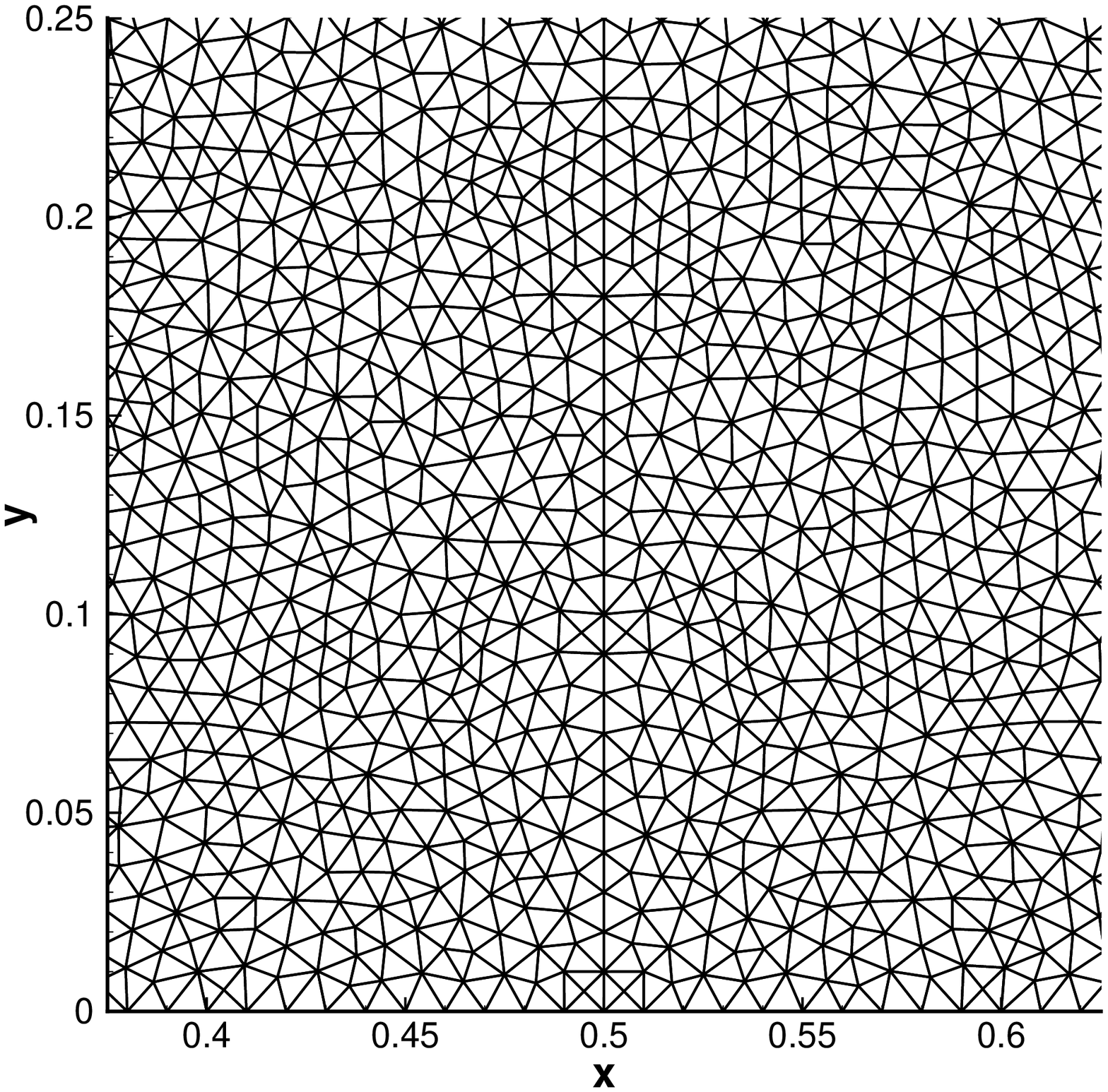}
\caption{\label{1d-shocktube-mesh2} One dimensional Riemann problem: triangular regular (left) and triangular irregular (right) meshes with.}
\end{figure}

\begin{figure}[!htb]
\centering
\includegraphics[width=0.45\textwidth]{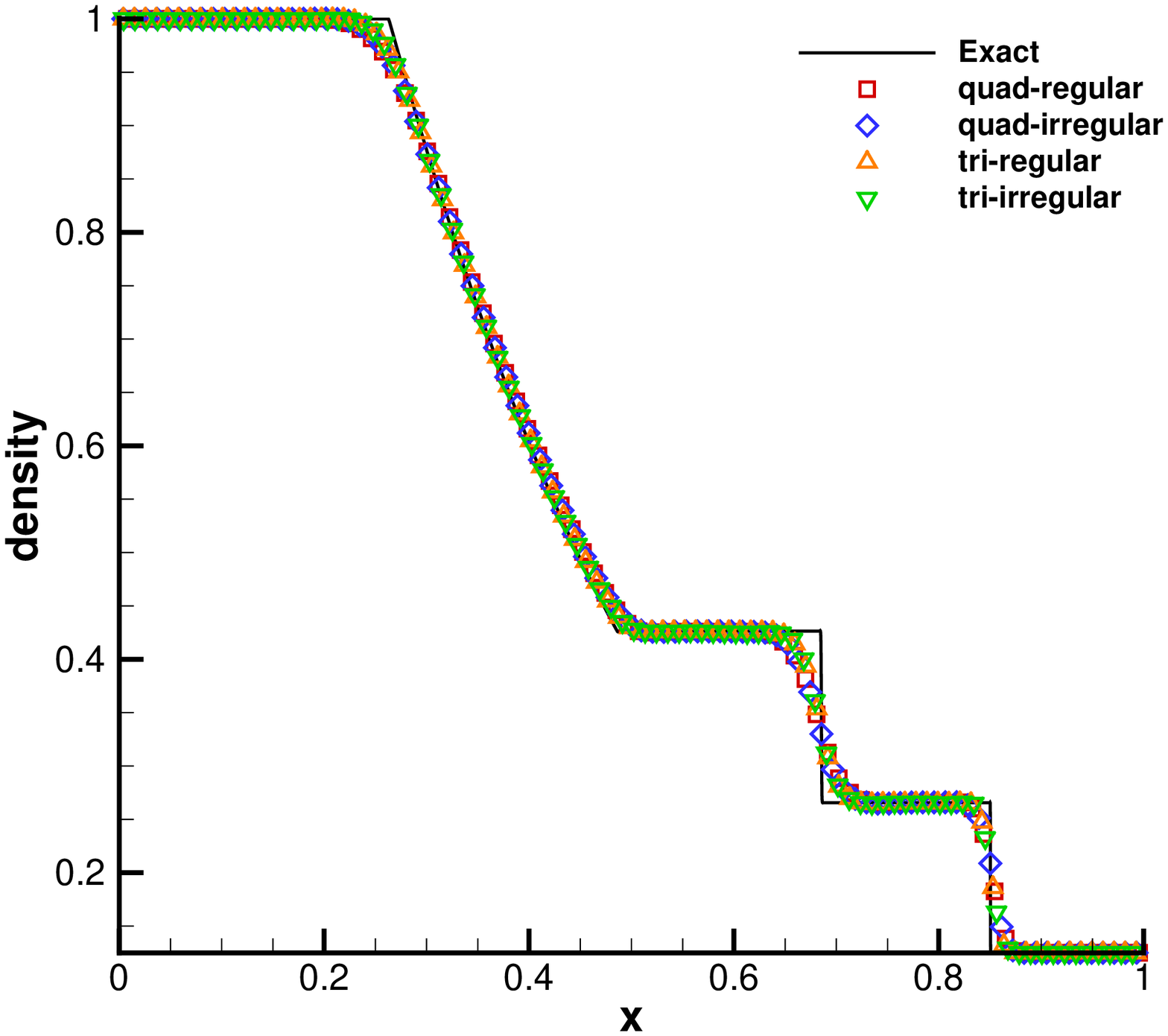}\includegraphics[width=0.45\textwidth]{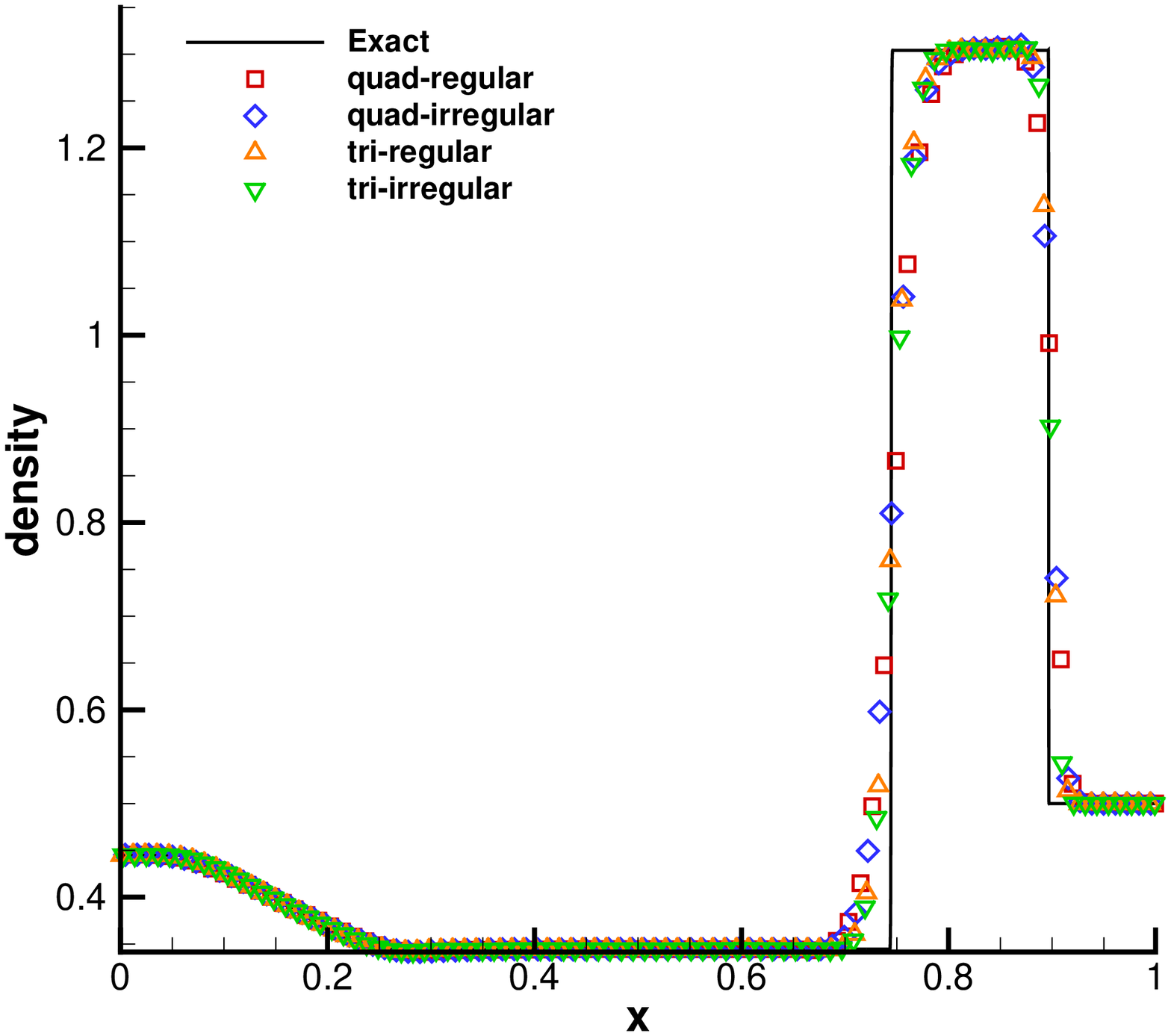}
\includegraphics[width=0.45\textwidth]{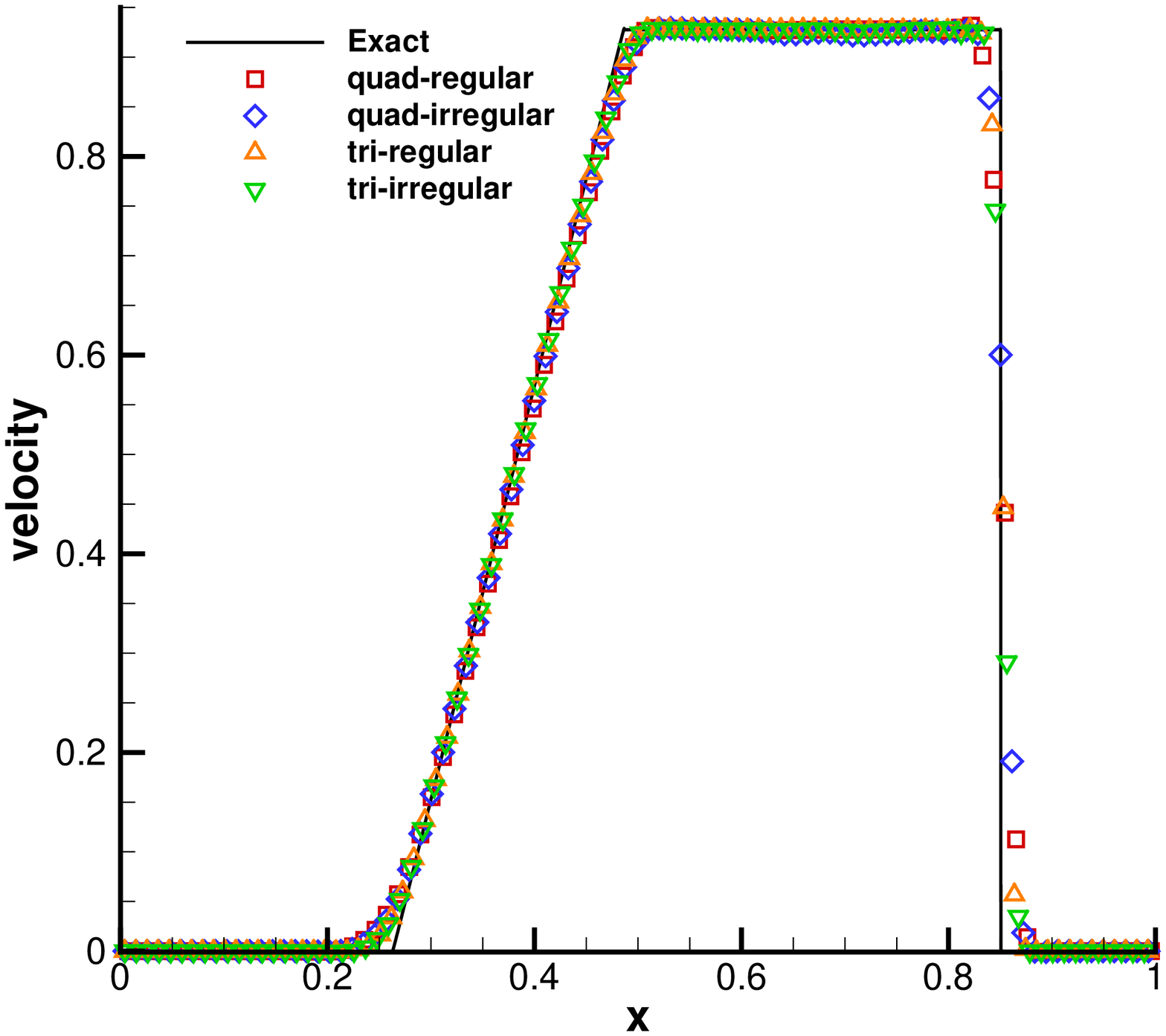}\includegraphics[width=0.45\textwidth]{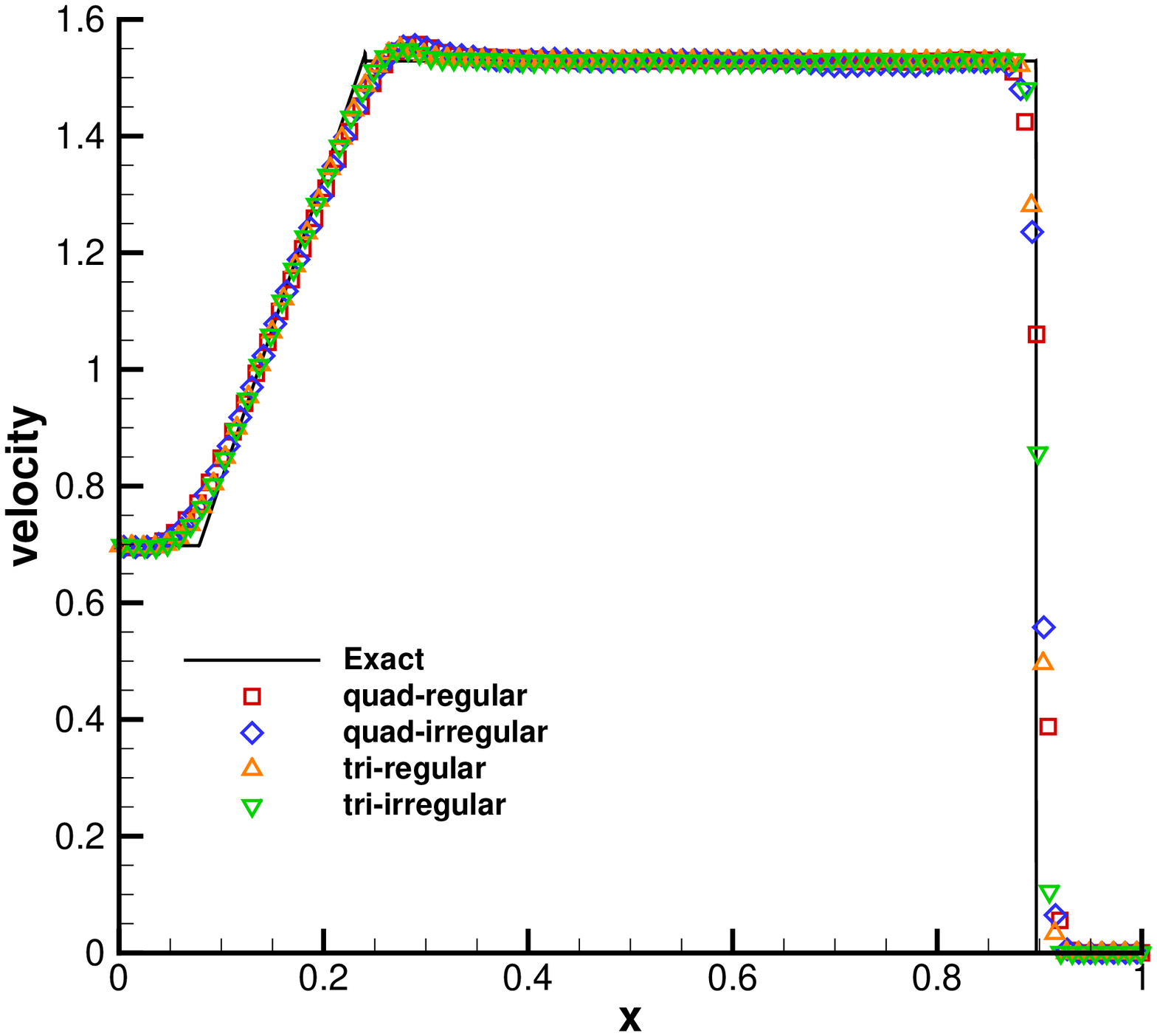}
\includegraphics[width=0.45\textwidth]{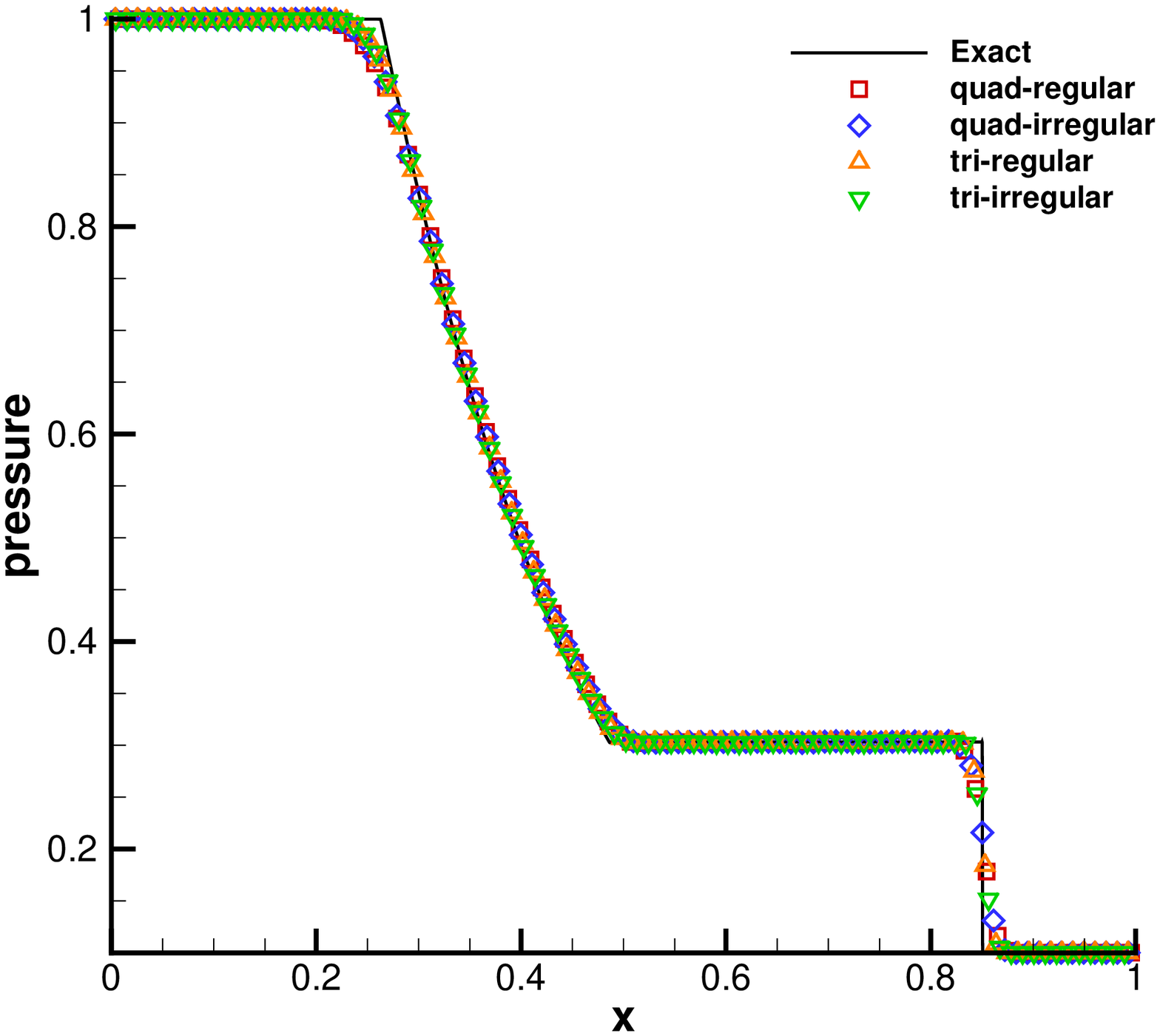}\includegraphics[width=0.45\textwidth]{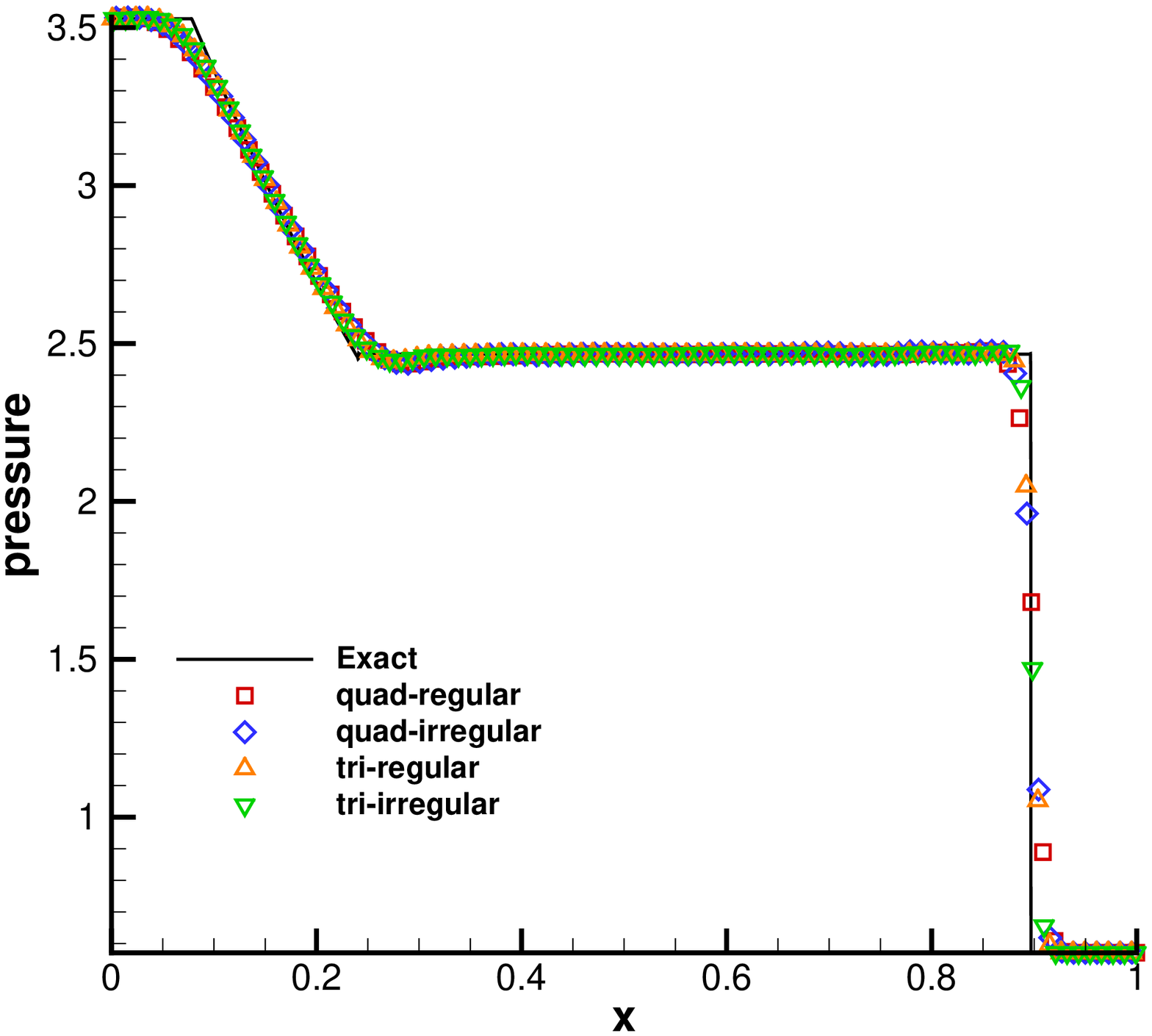}
\caption{\label{1d-shocktube-qua} One dimensional Riemann problem: the density, velocity and pressure distribution  for Sod problem at $t=0.2$ (left) and Lax problem at $t=0.16$ (right) on quadrilateral and triangular meshes.}
\end{figure}

\subsection{One dimensional Riemann problem}
In this case, four one-dimensional Riemann problems are tested by the
third-order WENO scheme on the quadrilateral and triangular meshes.
The first one is the Sod problem, and the initial condition is given
as follows
\begin{equation*}
(\rho,U,p) = \begin{cases}
(1,0,1),  0\leq x<0.5,\\
(0.125, 0, 0.1), 0.5\leq x<1.
\end{cases}
\end{equation*}
The second one is the Lax problem, and the initial condition is
given as follows
\begin{equation*}
(\rho,U,p) = \begin{cases}
(0.445, 0.698, 3.528),   0\leq x<0.5,\\
(0.5, 0, 0.571), 0.5\leq x<1.
\end{cases}
\end{equation*}
For these two cases, the computational domain is
$[0,1]\times[0,0.5]$, and the mesh size is $h=1/100$. Non-reflection
boundary condition is adopted at the boundaries of the computational
domain. The local quadrilateral and triangular meshes are given in
Figure.\ref{1d-shocktube-mesh1} and Figure.\ref{1d-shocktube-mesh2},
where both regular and irregular meshes are tested by the current
scheme. The numerical results for the Sod problem at $t=0.2$, and
for the Lax problem at $t=0.16$ are presented in
Figure.\ref{1d-shocktube-qua}. The density,
velocity and pressure distributions at the center horizontal line on
different meshes agree well with the exact solution. Due to more
number of cell with the same mesh scale, the scheme performs better
on triangular meshes than quadrilateral meshes.
 
\begin{figure}[!htb]
\centering
\includegraphics[width=0.485\textwidth]{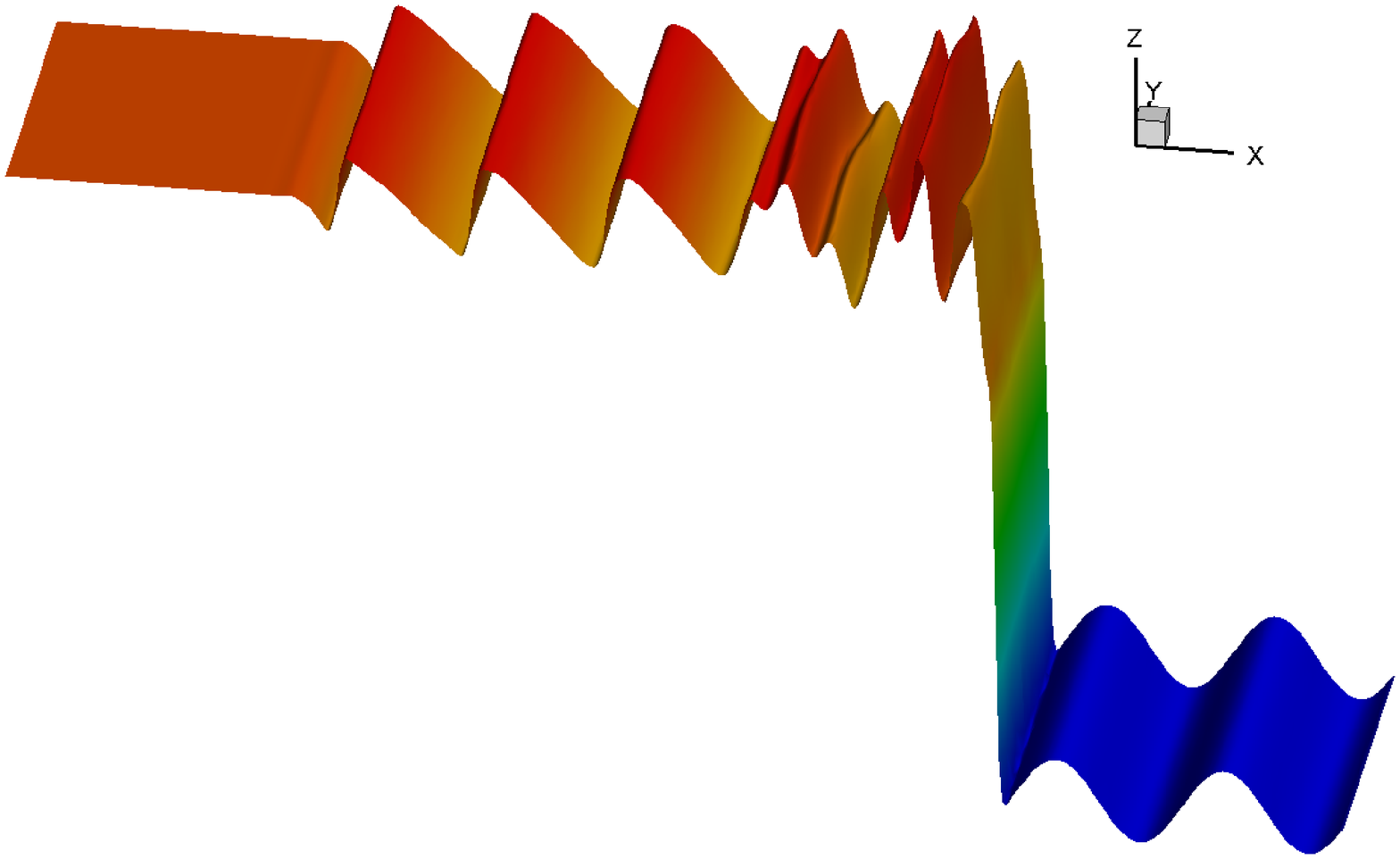}
\includegraphics[width=0.485\textwidth]{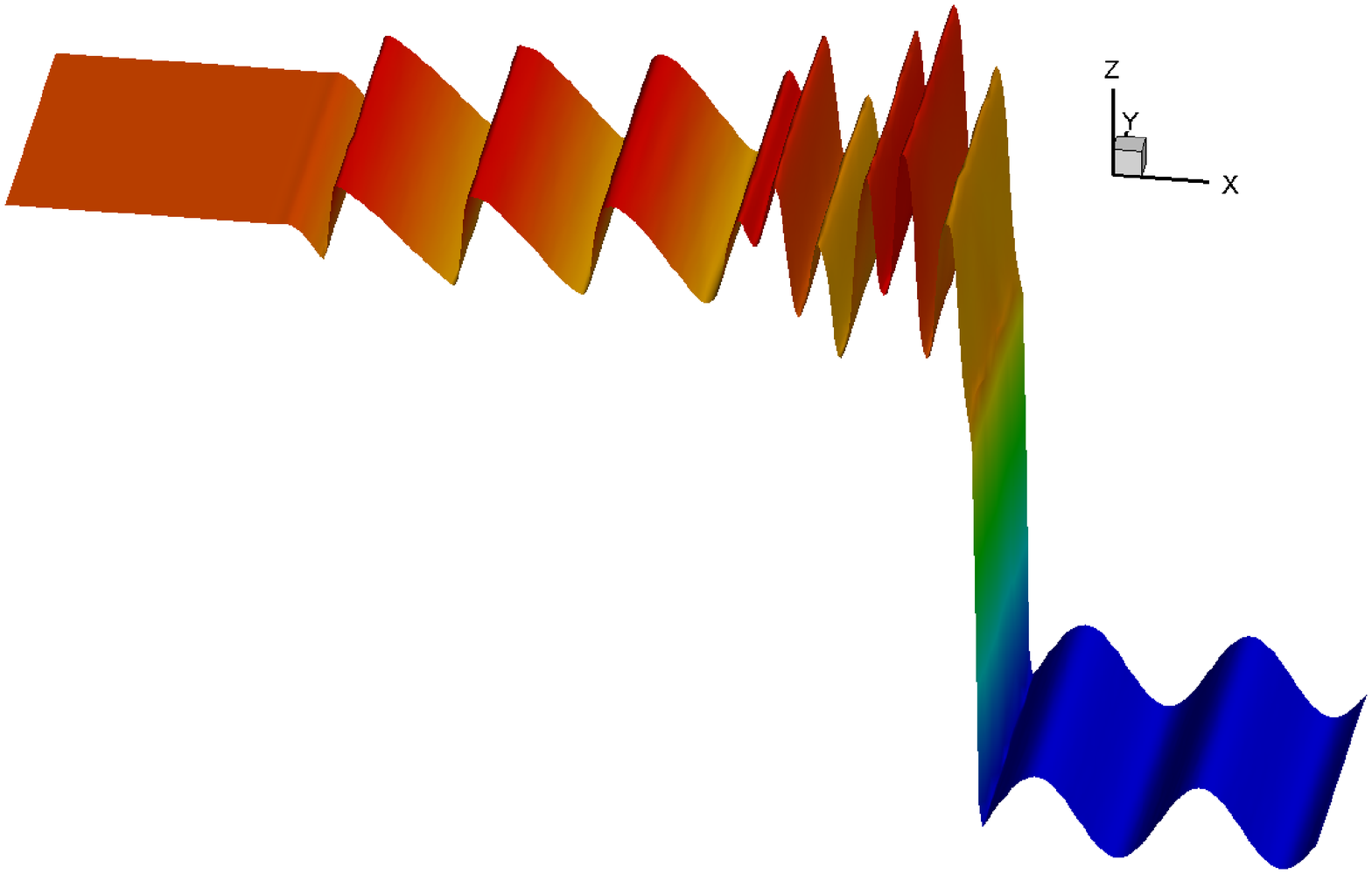}
\caption{\label{1d-shocktube-shu-1} One dimensional Riemann problem: the 3-D density distribution for Shu-Osher problem at $t=0.038$ on quadrilateral (left) and triangular (right) mesh with $h=1/40$.}
\includegraphics[width=0.485\textwidth]{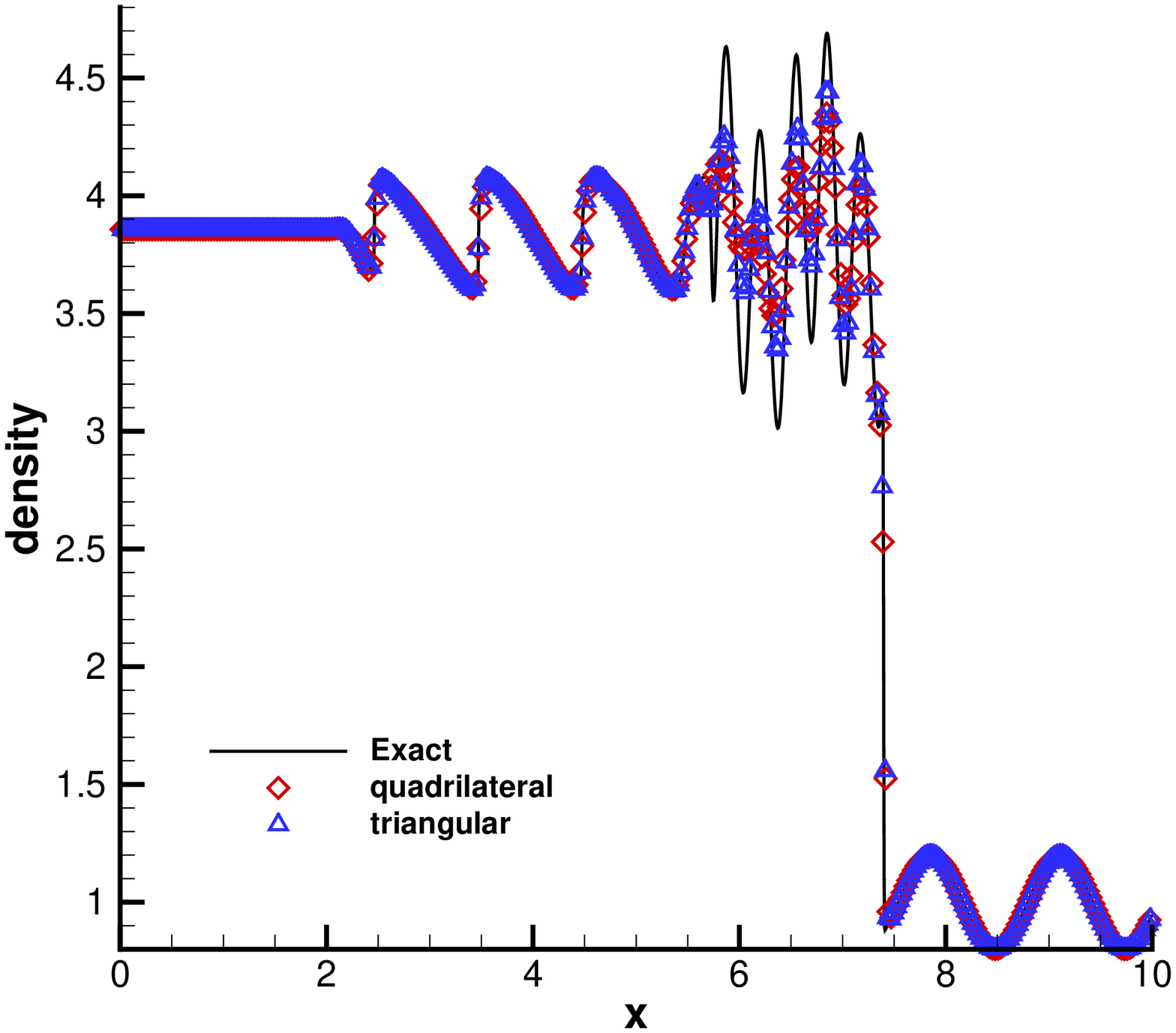}
\includegraphics[width=0.485\textwidth]{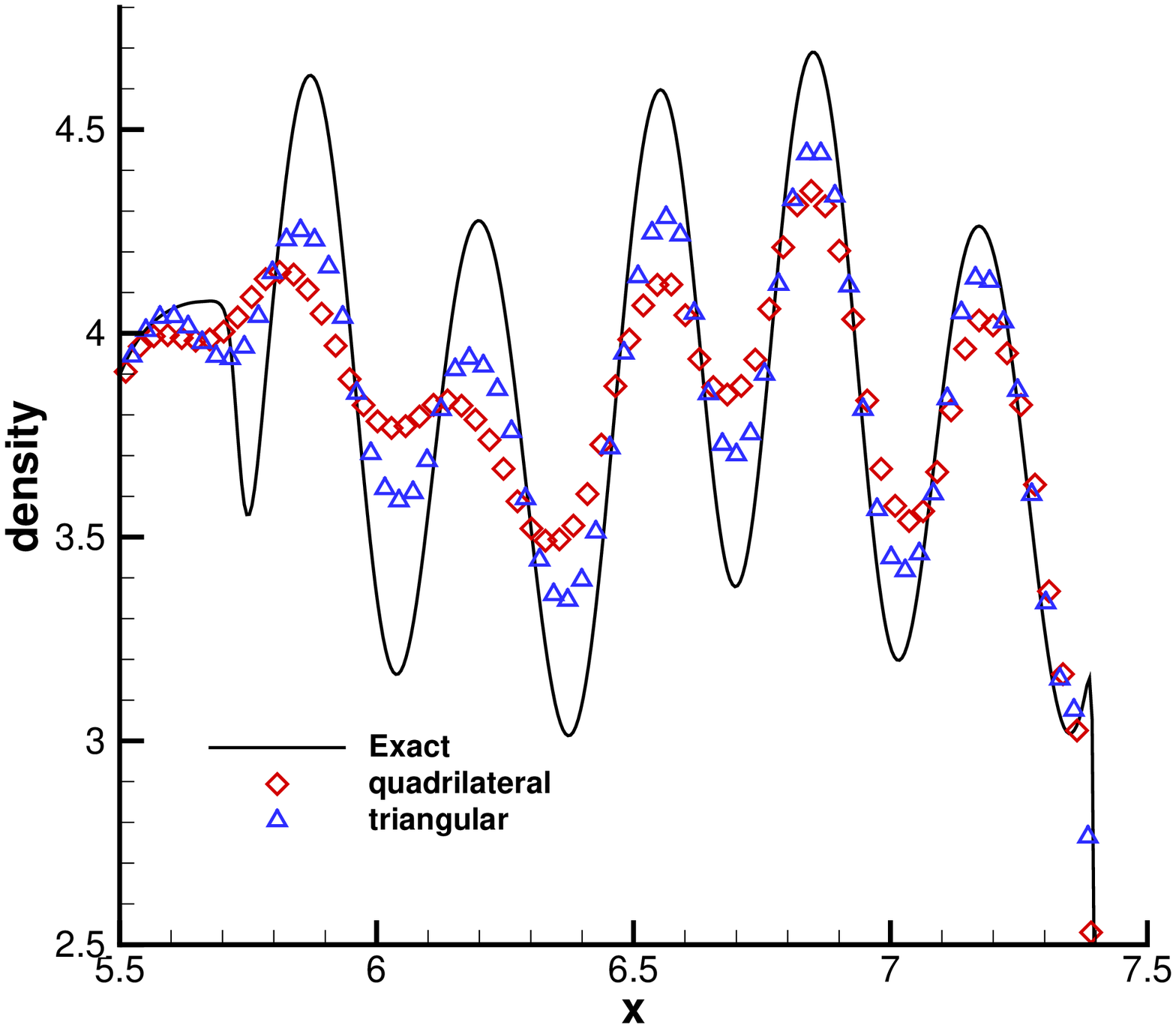}
\caption{\label{1d-shocktube-shu-2} One dimensional Riemann problem: the density distribution, local enlargement for Shu-Osher problem at $t=0.038$ on quadrilateral and triangular mesh with $h=1/40$.}
\end{figure}

To test the performance of capturing high frequency waves, the
Shu-Osher problem \cite{Case-Shu-Osher} is tested, which is a case
with the density wave interacting with shock. The initial condition
is given as follows
\begin{equation*}
(\rho,U,p)=\left\{\begin{array}{ll}
(3.857134, 2.629369, 10.33333),  \ \ \ \ &  x \leq 1,\\
(1 + 0.2\sin (5x), 0, 1),  &  1 <x.
\end{array} \right.
\end{equation*}
The computational domain is $[0, 10]\times[0, 0.25]$ and $h=1/40$
irregular quadrilateral and triangular meshes are used. The
reflected boundary condition is applied in the $y$ direction. The density distributions, local enlargement at the center line and 3-D distribution are presented in Figure.\ref{1d-shocktube-shu-1} and Figure.\ref{1d-shocktube-shu-2}  at
$t=1.8$. Due to the more cells of triangular mesh than quadrilateral
mesh with the same mesh scale, the result of triangular mesh
captures extremum better.

\begin{figure}[!htb]
\centering
\includegraphics[width=0.475\textwidth]{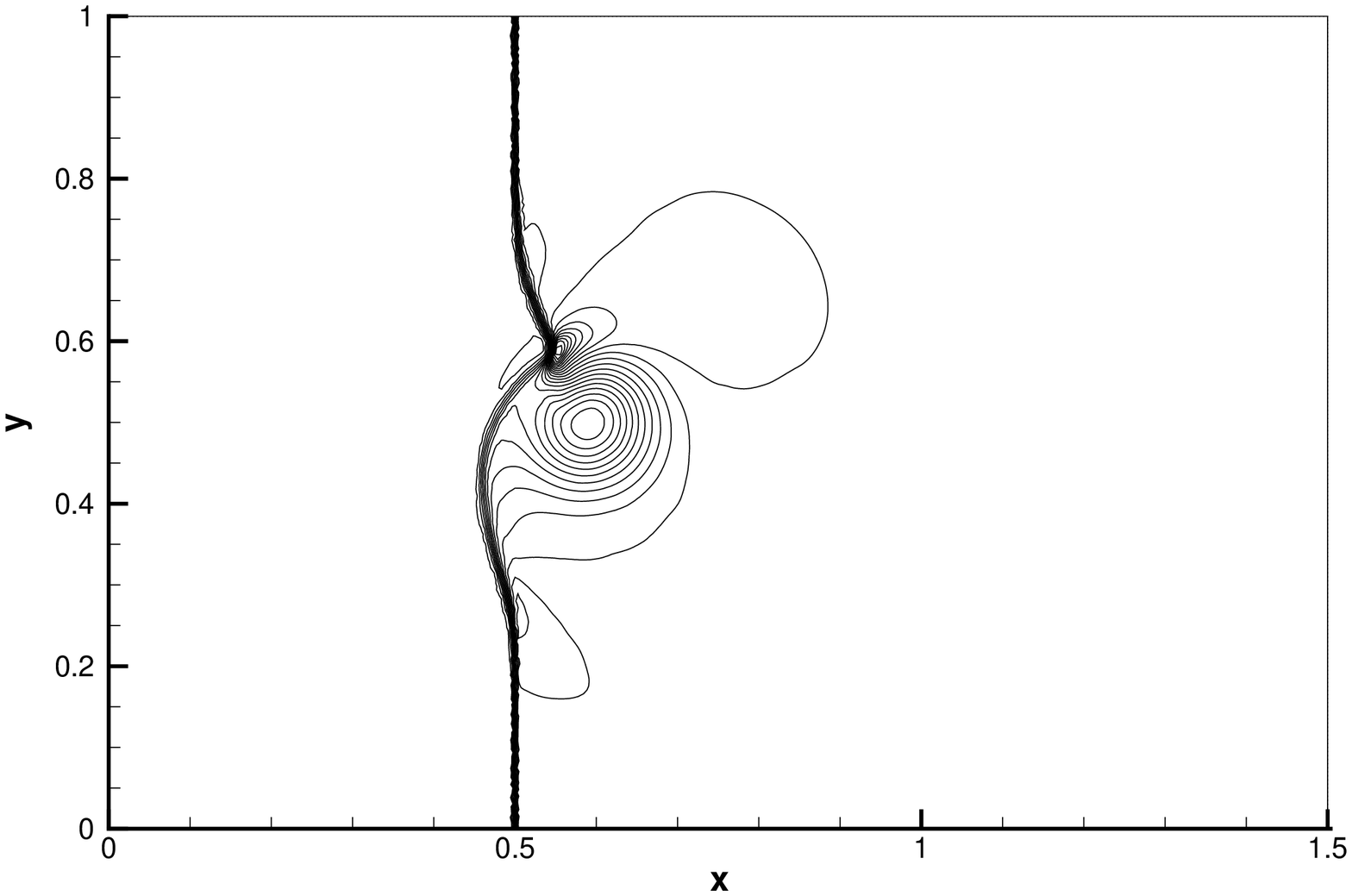}
\includegraphics[width=0.475\textwidth]{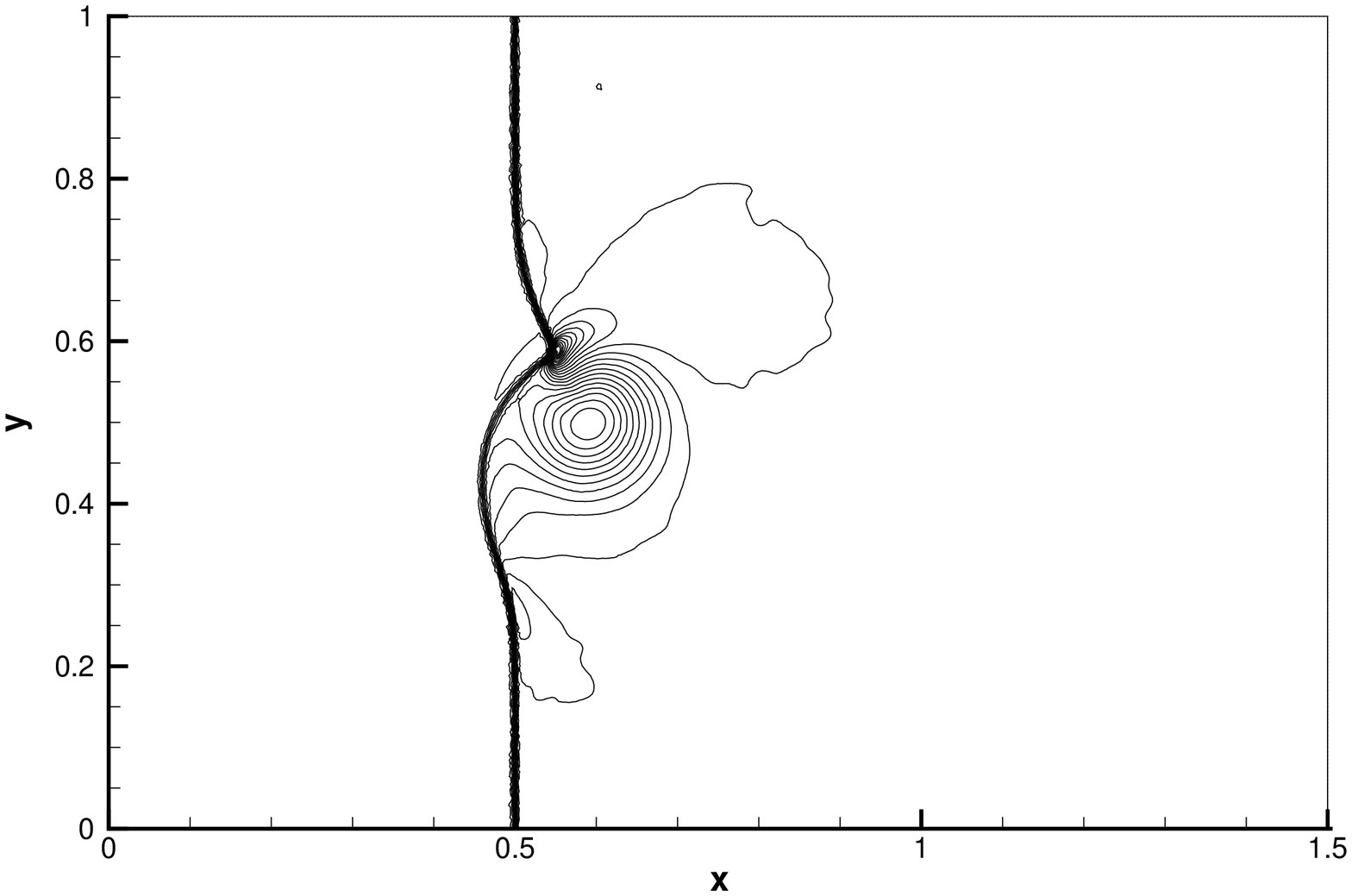}
\includegraphics[width=0.475\textwidth]{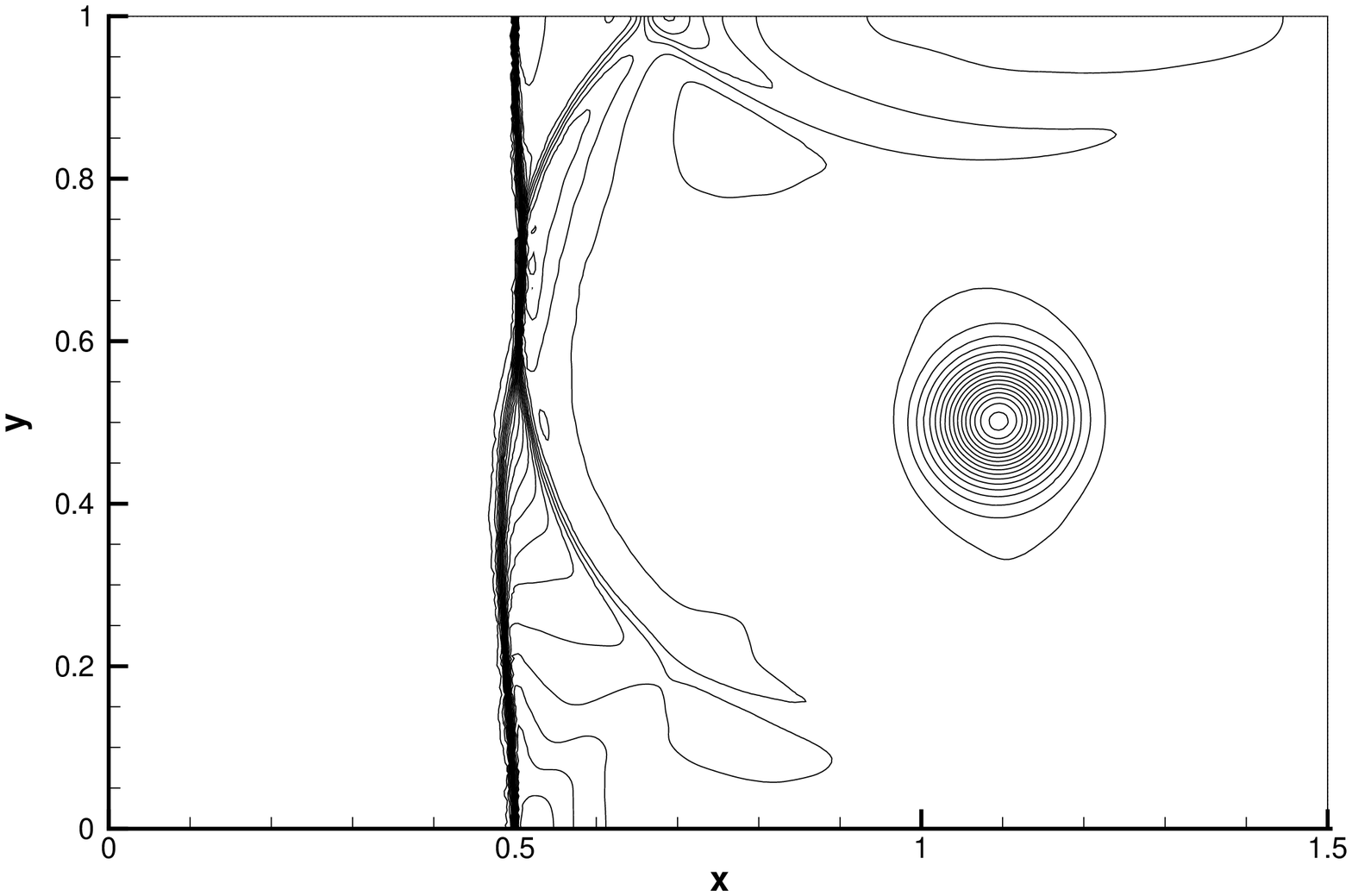}
\includegraphics[width=0.475\textwidth]{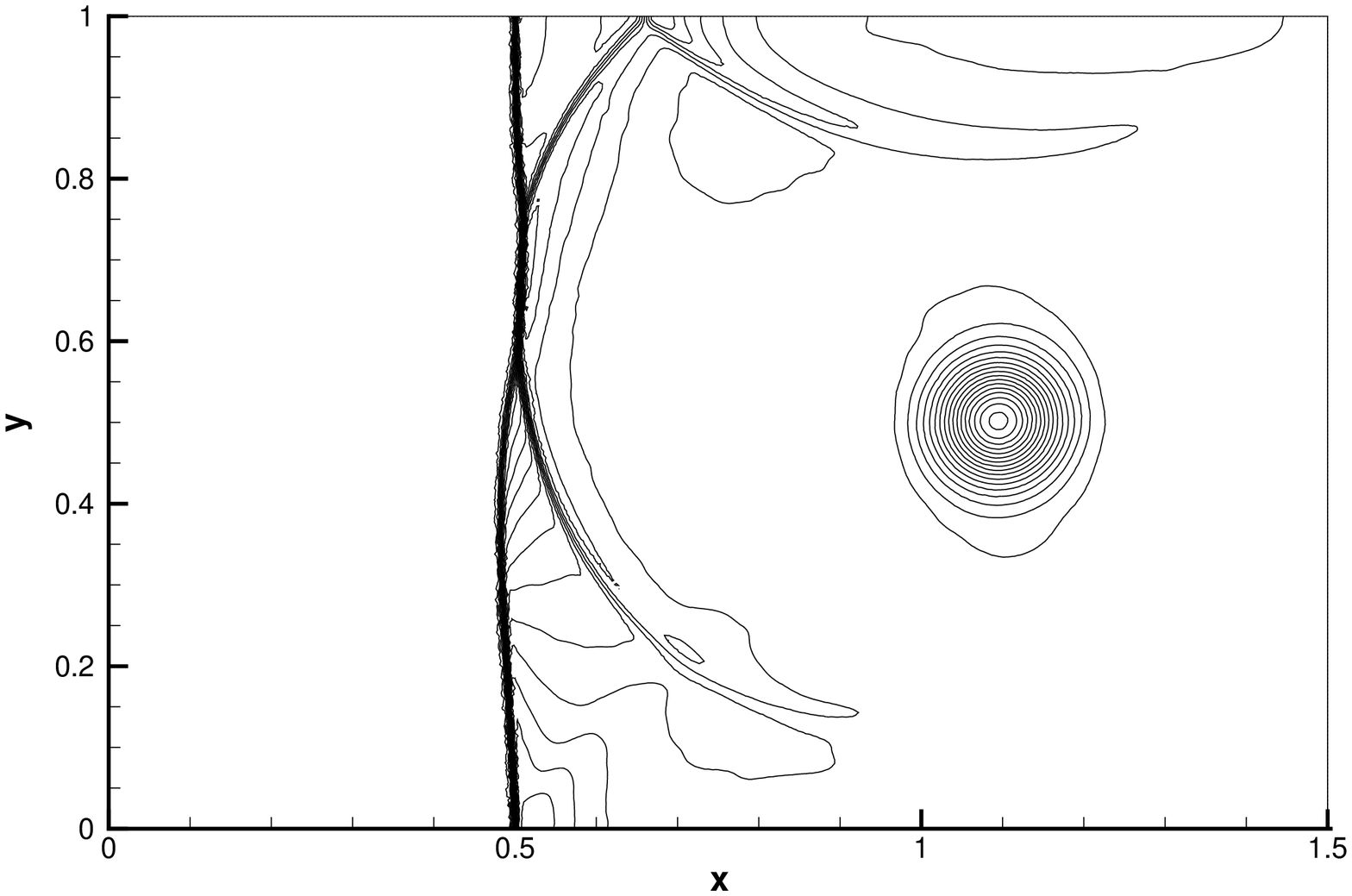}
\caption{\label{2d-Vorshock}Shock-vortex interaction: density
distribution on quadrilateral irregular mesh (left) and triangular irregular mesh
(right) at $t=0.3, 0.8$ with $h=1/200$.}
\end{figure}

\begin{figure}[!htb]
\centering
\includegraphics[width=0.49\textwidth]{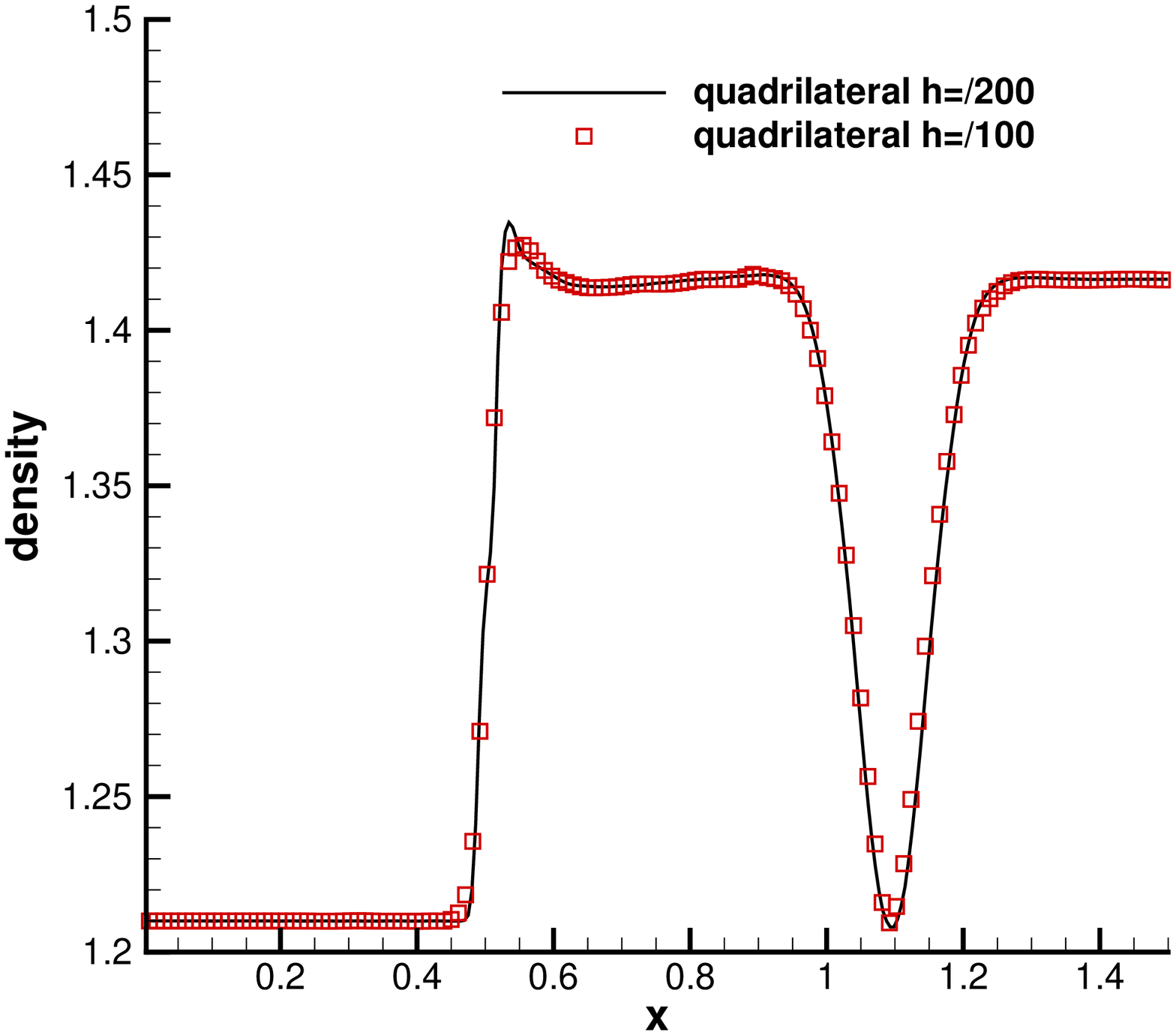}
\includegraphics[width=0.49\textwidth]{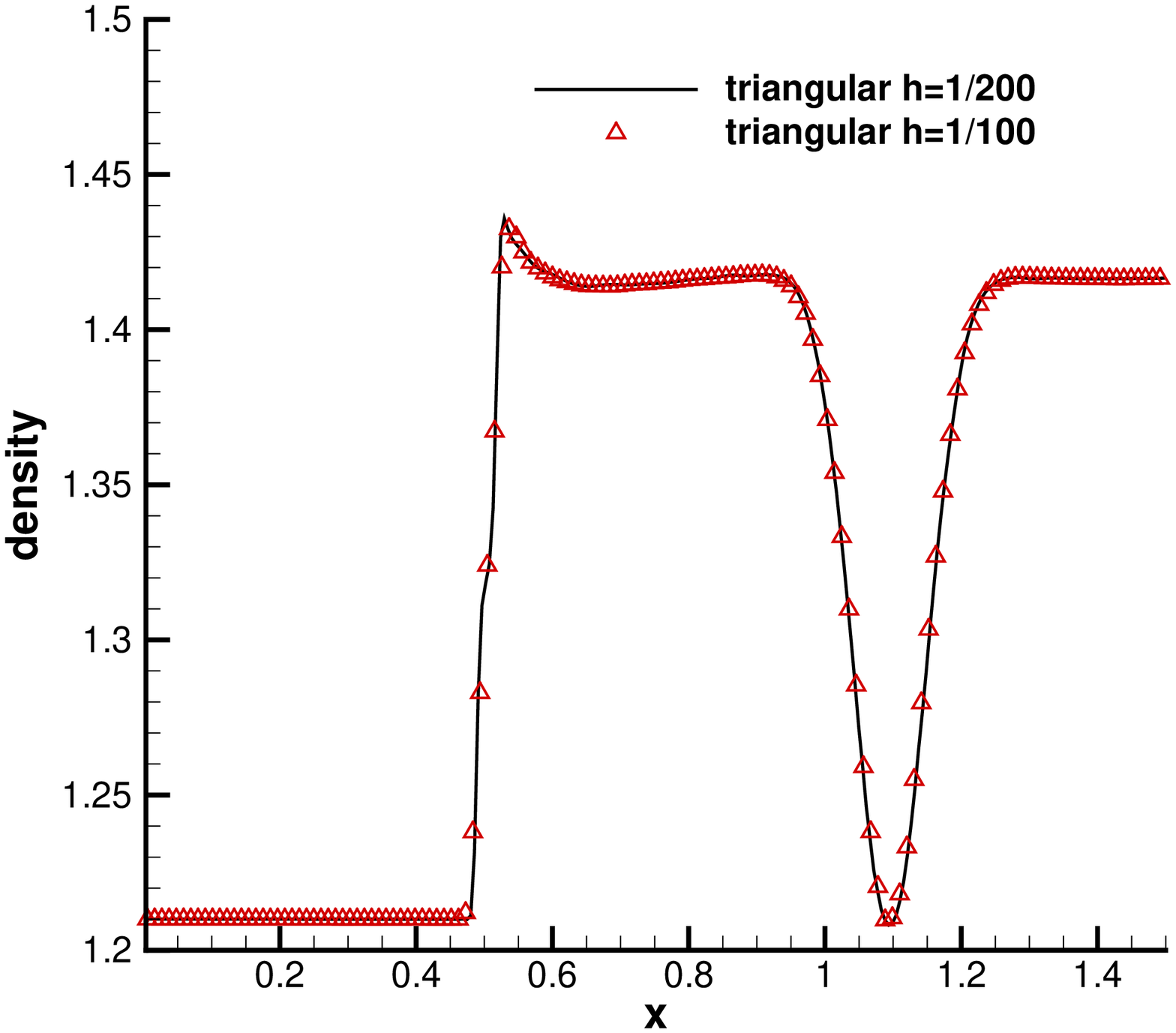}
\caption{\label{2d-Vorshock-pressure}Shock-vortex interaction:
density distribution on quadrilateral irregular meshes (left) and triangular irregular
meshes (right) with $h=1/100, 1/200$ at $t=0.8$.}
\end{figure}

\begin{figure}[!htb]
\centering
\includegraphics[width=0.495\textwidth]{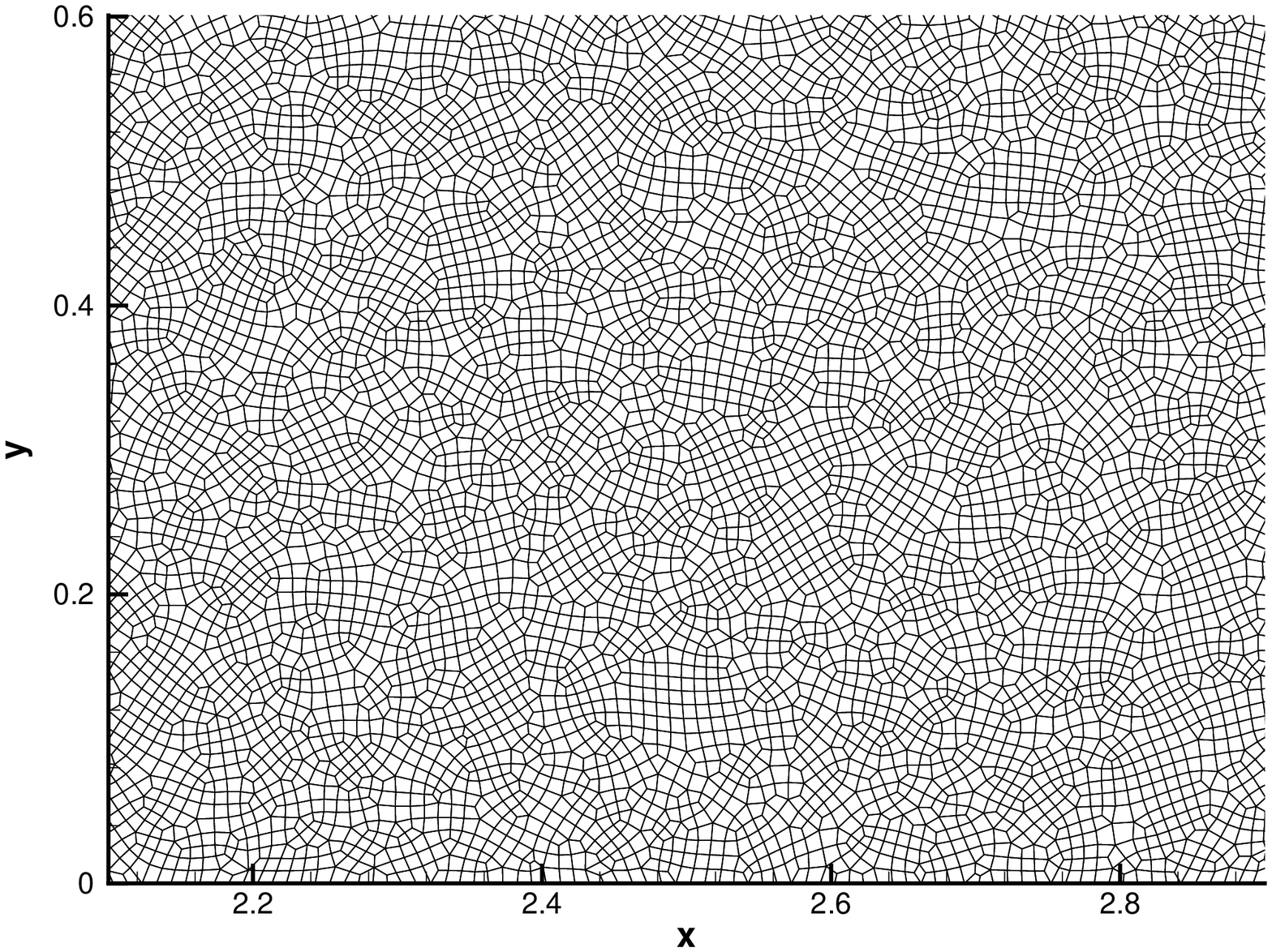}
\includegraphics[width=0.495\textwidth]{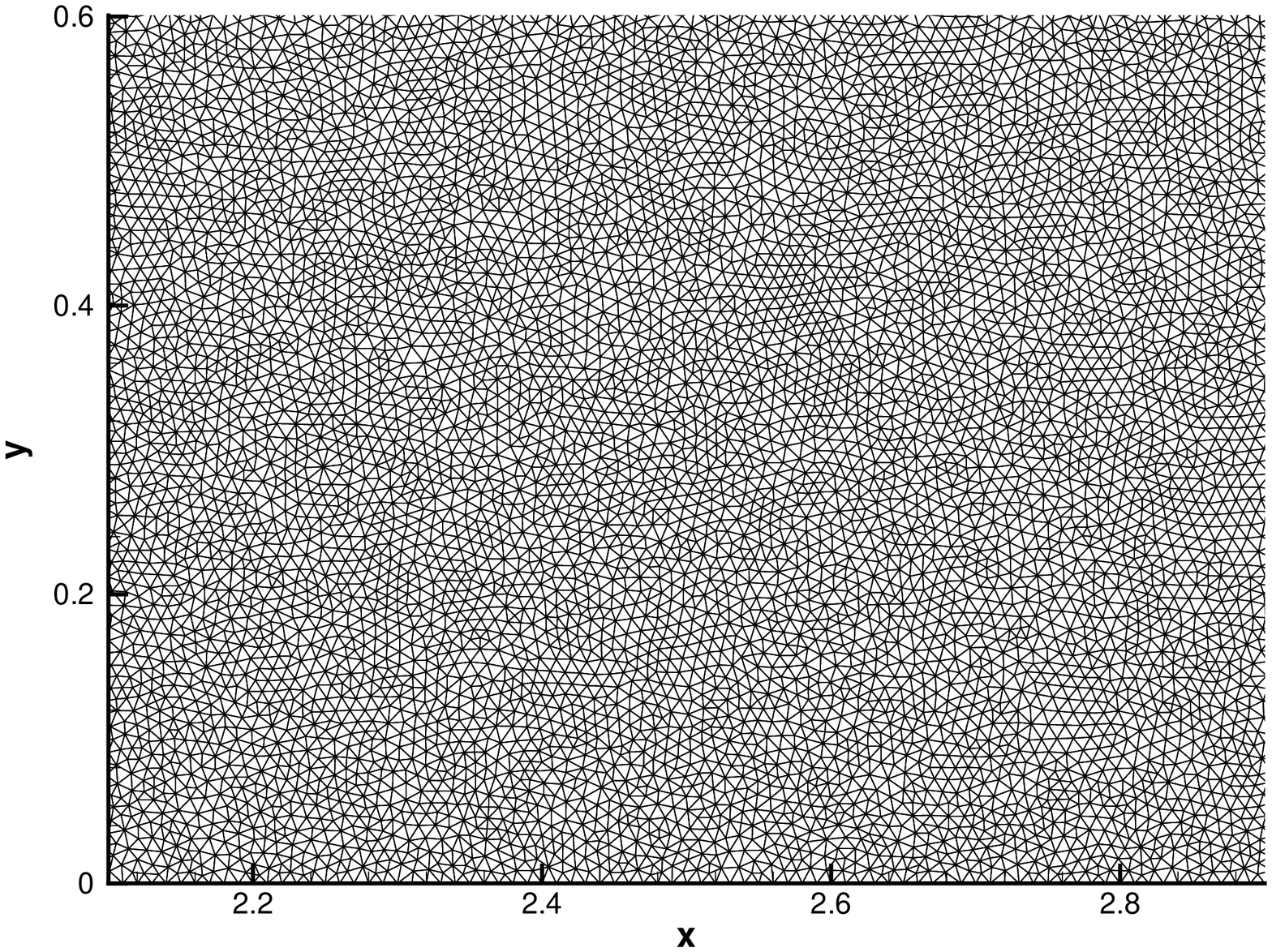}\\
\includegraphics[width=0.495\textwidth]{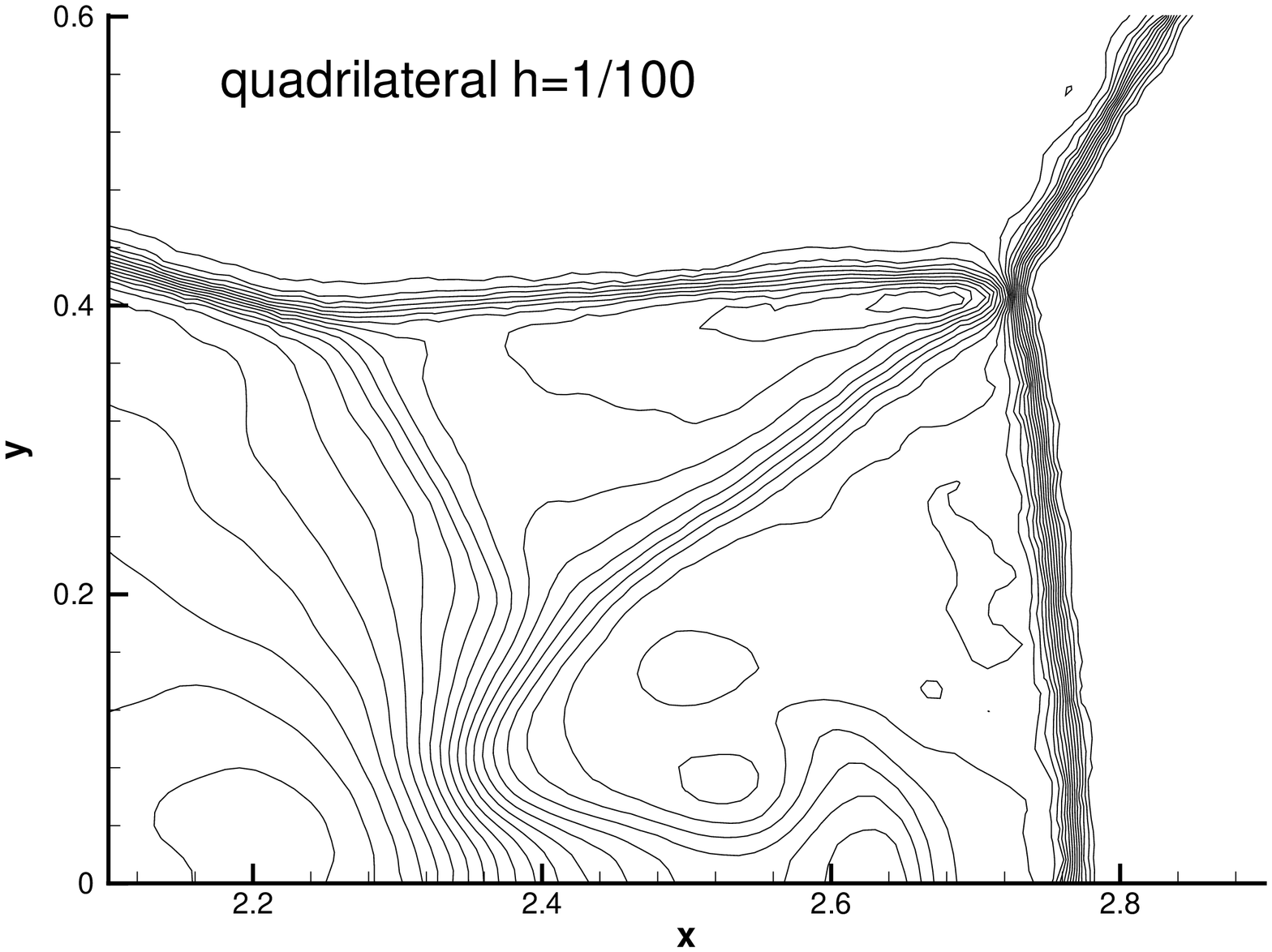}
\includegraphics[width=0.495\textwidth]{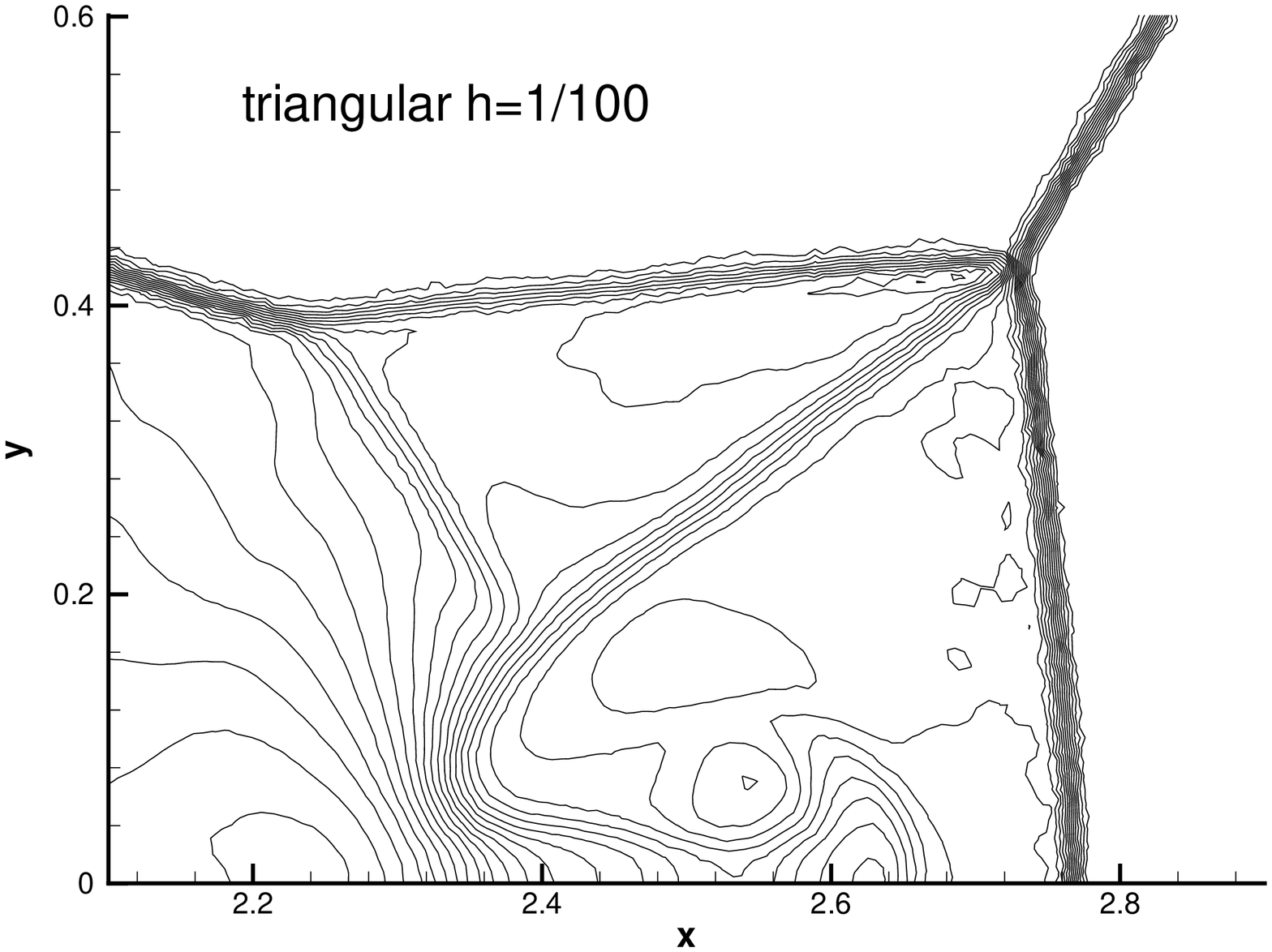}
\caption{\label{double-mach-1} Double Mach reflection: local density
and mesh distributions on quadrilateral mesh (left) and triangular
mesh (right) with $h=1/100$ at $t=0.2$.}
\end{figure}

\subsection{Shock vortex interaction}
The interaction between a stationary shock and a vortex for the
inviscid flow \cite{WENO-JS} is presented. The computational domain
is taken to be $[0, 2]\times[0, 1]$. A stationary Mach $1.1$ shock
is positioned at $x=0.5$ and normal to the $x$-axis. The left
upstream state is $(\rho, U, V, p) = (Ma^2,\sqrt{\gamma}, 0, 1)$,
where $Ma$ is the Mach number. A small vortex is obtained through a
perturbation on the mean flow with the velocity $(U, V)$,
temperature $T=p/\rho$, and entropy $S=\ln(p/\rho^\gamma)$.  The
perturbation is expressed as
\begin{align*}
&(\delta U,\delta V)=\kappa\eta e^{\mu(1-\eta^2)}(\sin\theta,-\cos\theta),\\
&\delta
T=-\frac{(\gamma-1)\kappa^2}{4\mu\gamma}e^{2\mu(1-\eta^2)},\delta
S=0,
\end{align*}
where $\eta=r/r_c$, $r=\sqrt{(x-x_c)^2+(y-y_c)^2}$, and $(x_c,
y_c)=(0.25, 0.5)$ is the center of the vortex. Here $\kappa$
indicates the strength of the vortex, $\mu$ controls the decay rate
of the vortex, and $r_c$ is the critical radius for which the vortex
has the maximum strength. In the computation, $\kappa=0.3$,
$\mu=0.204$, and $r_c=0.05$. The reflected boundary conditions are
used on the top and bottom boundaries. The density distributions
with mesh size $h=1/200$ at $t=0.3$ and $0.8$ are shown in
Figure.\ref{2d-Vorshock}. By $t=0.8$, one branch of the shock
bifurcations has reached the top boundary and been reflected, and
the reflection is well captured. The detailed density distributions
along the center horizontal line on quadrilateral and triangular
meshes with $h=1/100$ and $1/200$ at $t=0.8$ are shown in
Figure.\ref{2d-Vorshock-pressure}.

\begin{figure}[!htb]
\centering
\includegraphics[width=0.495\textwidth]{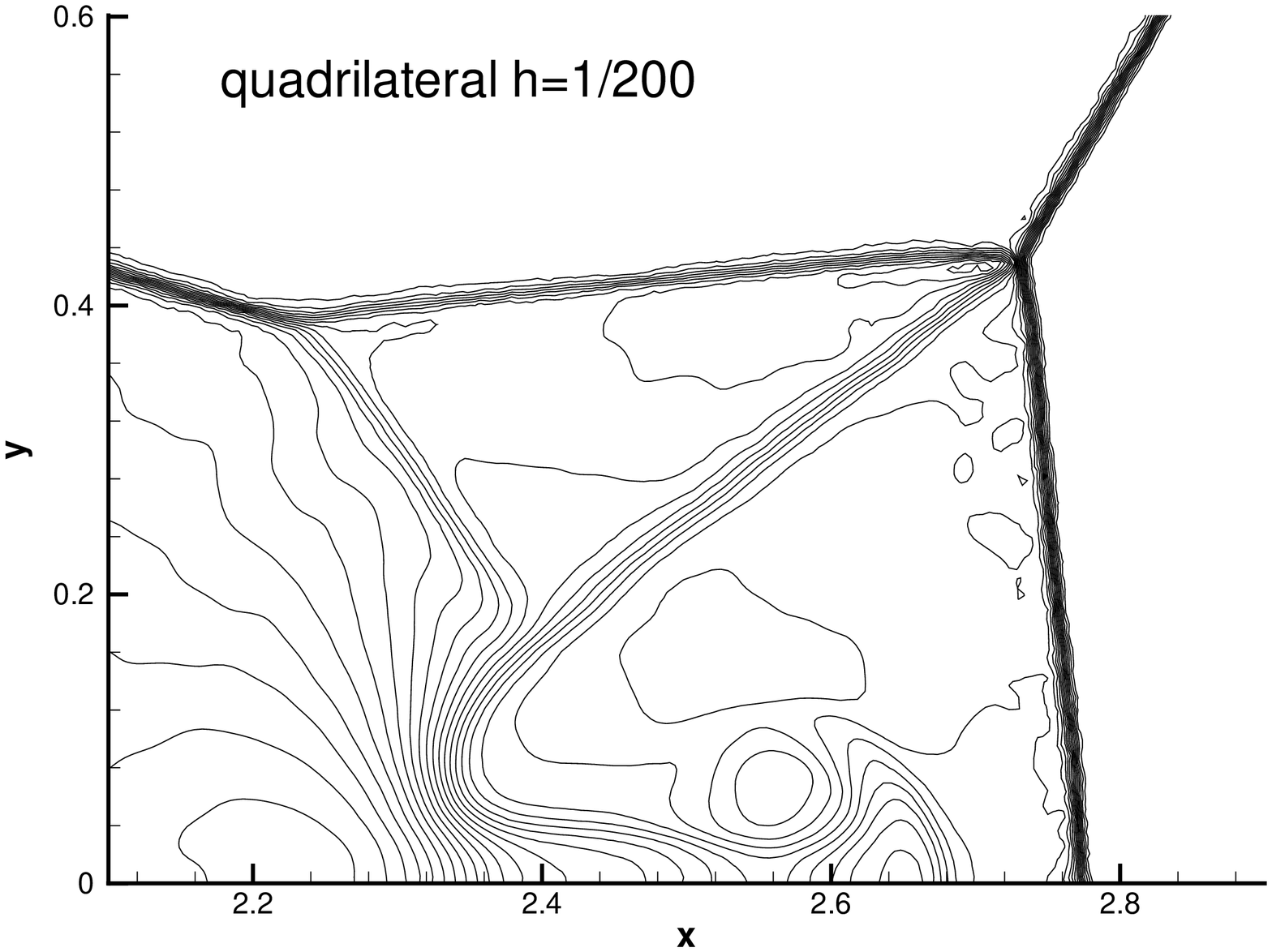}
\includegraphics[width=0.495\textwidth]{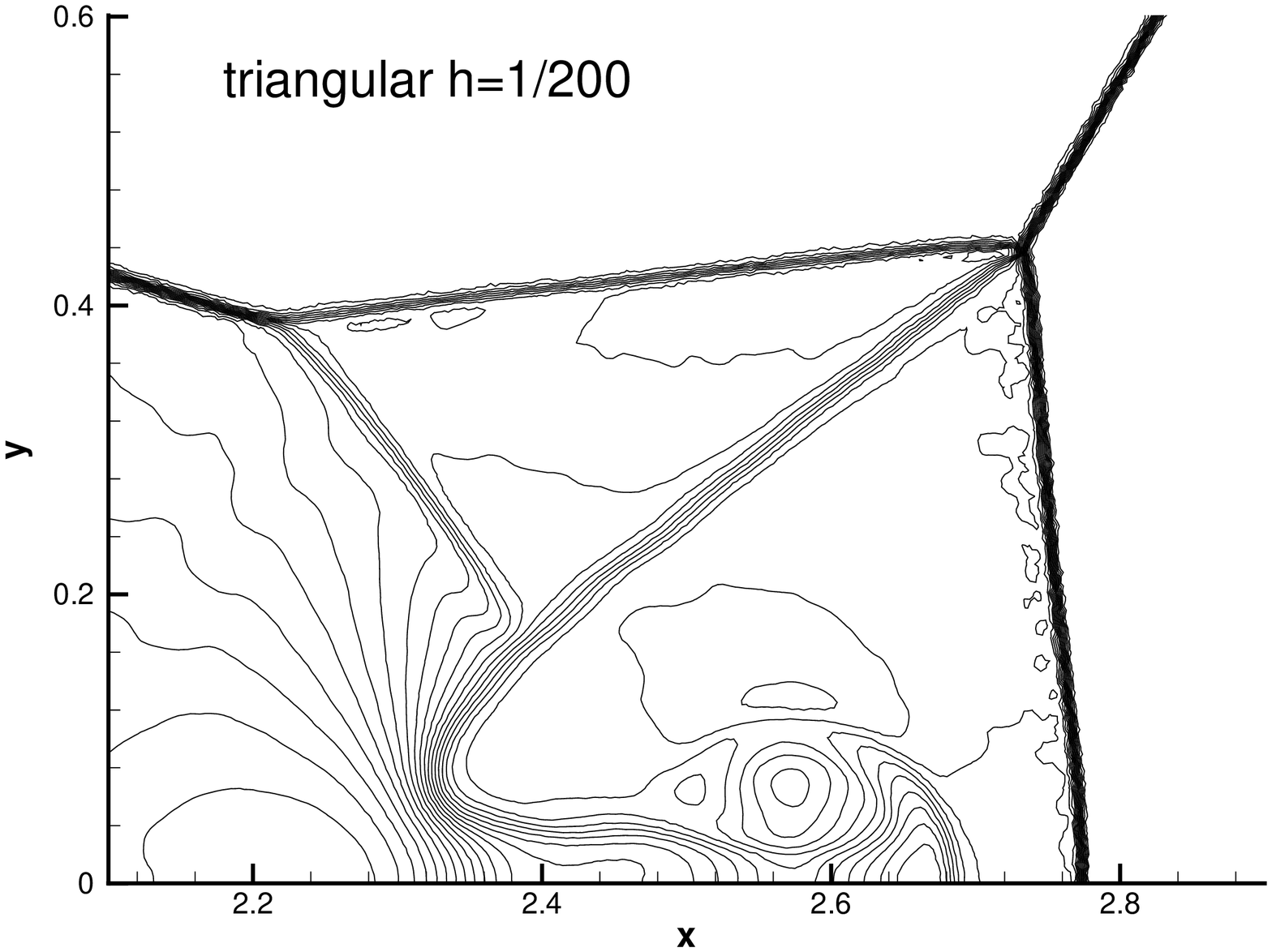}\\
\includegraphics[width=0.495\textwidth]{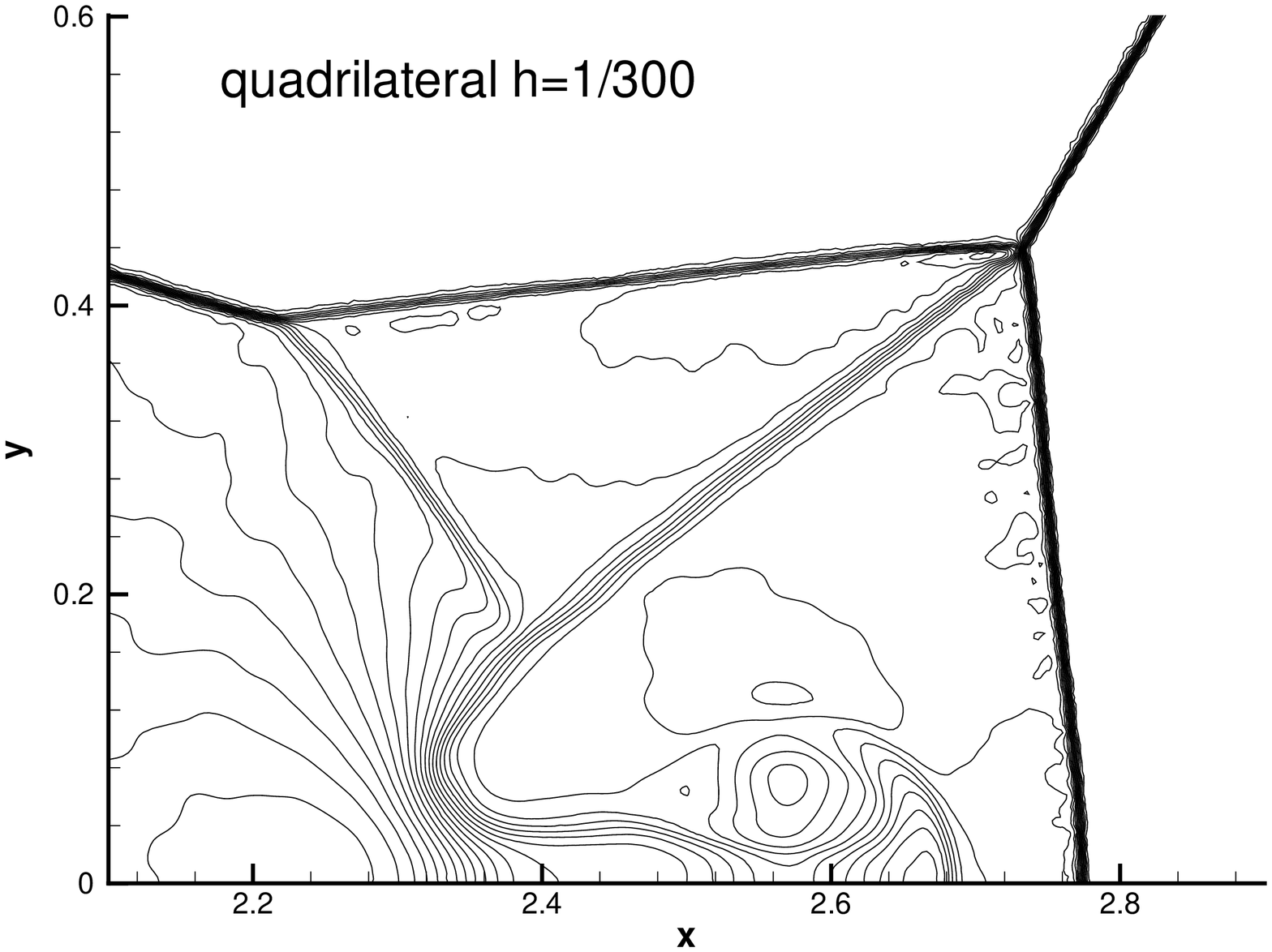}
\includegraphics[width=0.495\textwidth]{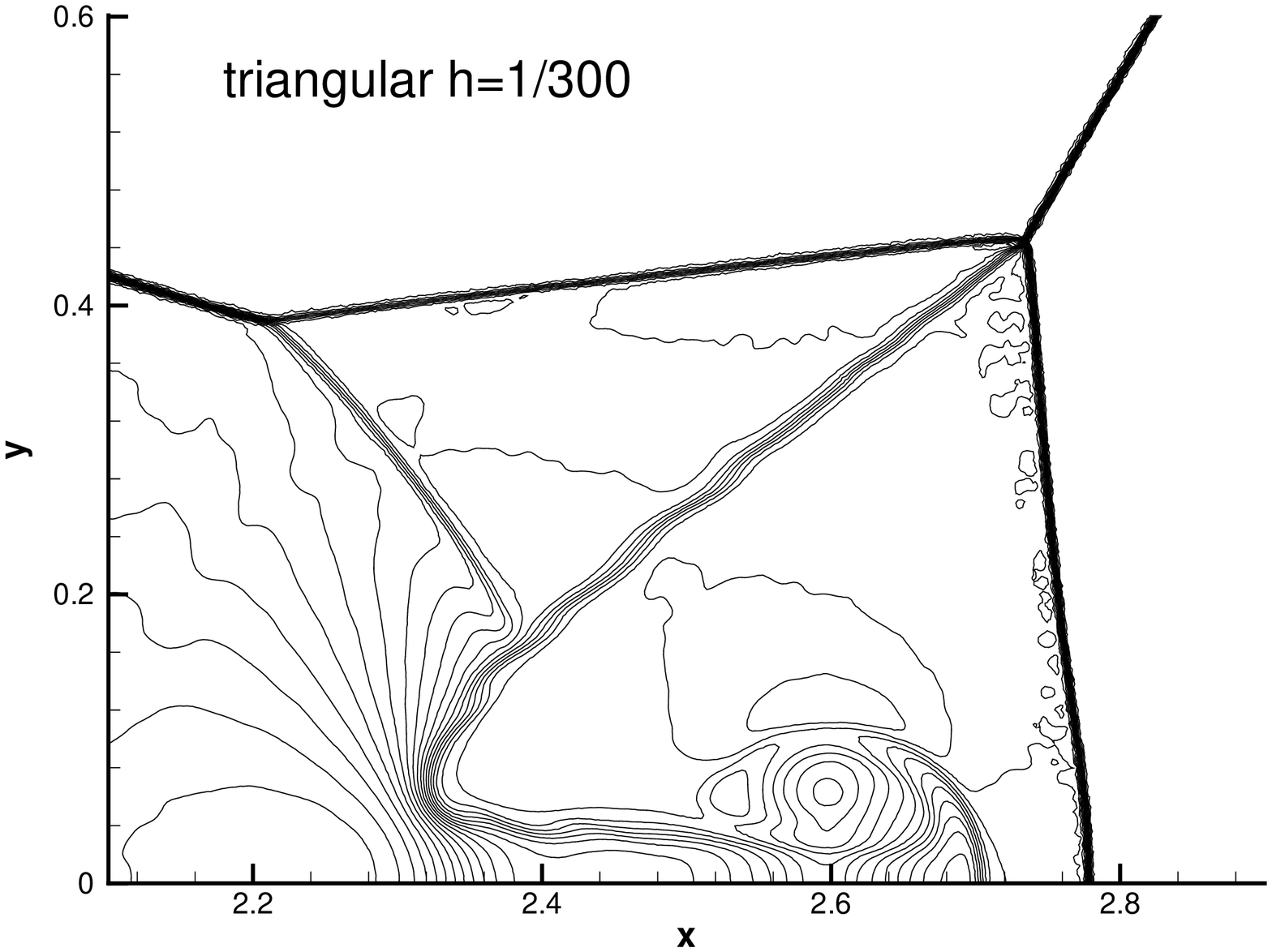}
\caption{\label{double-mach-2} Double Mach reflection: density
distributions on quadrilateral meshes (left) and triangular meshes
(right) with $h=1/200$ and $1/400$ at $t=0.2$.}
\end{figure}

\subsection{Double Mach reflection problem}
This problem was extensively studied by Woodward and Colella
\cite{Woodward-Colella} for the inviscid flow. The computational
domain is $[0,4]\times[0,1]$, and a solid wall lies at the bottom of
the computational domain starting from $x =1/6$. Initially a
right-moving Mach 10 shock is positioned at $(x,y)=(1/6, 0)$, and
makes a $60^\circ$ angle with the x-axis. The initial pre-shock and
post-shock conditions are
\begin{align*}
(\rho, U, V, p)&=(8, 4.125\sqrt{3}, -4.125,
116.5),\\
(\rho, U, V, p)&=(1.4, 0, 0, 1).
\end{align*}
The reflective boundary condition is used at the wall, while for the
rest of bottom boundary, the exact post-shock condition is imposed.
At the top boundary, the flow values are set to describe the exact
motion of the Mach $10$ shock. For this case, the current WENO
scheme is tested by both quadrilateral and triangular meshes with
mesh size $h=1/100, 1/200$ and $1/300$, and the local meshes near
the Mach stem is given in Fig.\ref{double-mach-1} with $h=1/100$.
The density distributions with mesh size $h=1/100, 1/200$ and
$1/400$ at $t=0.2$ are presented in Figure.\ref{double-mach-1} and
Figure.\ref{double-mach-2}, respectively. The current scheme resolves
the flow structure under the triple Mach stem clearly with the mesh
refinement, and the robustness of the current scheme is well
validated.

\section{Conclusion}
In this paper, a new third-order WENO scheme has been developed
based on both unstructured quadrilateral and triangular meshes for
the hyperbolic conservation laws. As a starting point of WENO
reconstruction, a unified linear scheme was constructed for both
quadrilateral and triangular meshes with any local topology.
However, the very large linear weights may appear for the mesh with
lower quality, which will make the scheme unstable even for the
smooth flows. In the current reconstruction, an optimization
approach was introduced to cure the very large linear weights on
unstructured meshes. The splitting technique is considered to deal
with the negative weights obtained by the optimization approach.
With the optimization approach for very large weights and the
splitting technique for negative weights, the scheme becomes more
robust than the classical WENO scheme \cite{WENO-JS,Hu-Shu} with the
low quality meshes. The accuracy and expected convergence order,
which are independent of mesh quality, are obtained for both
quadrilateral and triangular meshes, and the stability of the
current scheme is not affected by the mesh quality.  A variety of
numerical tests with strong discontinuities are presented to
validate the robustness of the current scheme.

\section*{Acknowledgements}
The work of L. Pan is supported by NSFC (11701038) and China Postdoctoral Science
Foundation (2016M600065).  The work of S.H. Wang is supported by
NSAF (U1630247) and  NSFC (915303108).

\end{document}